\documentstyle[11pt]{article}
\input epsf

\newtheorem{example}{Example}[section]
\newtheorem{note}[example]{Note}
\newtheorem{theorem}[example]{Theorem}
\newtheorem{corollary}[example]{Corollary}
\newtheorem{definition}[example]{Definition}
\newtheorem{proposition}[example]{Proposition}

\newtheorem{lemma}[example]{Lemma}

\def\La{\Lambda}

\def\P{{\cal P}}

\def\wt{\mbox{\sl wt}\,}

\def\Z{{\bf Z}}

\def\blambda{{\hbox{\boldmath$\lambda$}}}
\def\bmu{{\hbox{\boldmath$\mu$}}}

\def\btau{{\hbox{\boldmath$\tau$}}}
\def\sblambda{{\hbox{\boldmath$\lambda$}}}
\def\sbmu{{\hbox{\boldmath$\mu$}}}

\def\N{{\ssym N}}

\def\Z{{\ssym Z}}
\def\Q{{\ssym Q}}
\def\C{{\ssym C}}
\def\JS{{\rm JS}}
\def\Proof{\medskip\noindent {\it Proof: }}
\def\cqfd{\hfill $\Box$ \bigskip}
\def\endex{\hfill $\diamond$ \medskip}
\def\adots{\mathinner{\mkern2mu\raise1pt\hbox{.} 
\mkern3mu\raise4pt\hbox{.}\mkern1mu\raise7pt\hbox{.}}}
\def\<{\langle}
\def\>{\rangle}
\def\resp{{\it resp.$\ $}}

\def\SG{{\goth S}}

\def\mod{{\rm\,mod\,}}

\def\LP{\Lambda^\prime}
\def\LPP{\Lambda^{\prime\prime}}

\def\Sl{{\goth sl}}
\def\glchap{\widehat{\hbox{\goth gl}}}
\def\slchap{\widehat{\goth sl}}
\def\h{{\goth h}}

\def\v{\mbox{\boldmath$v$}}
\def\F{{\cal F}}

\def\H{{\cal H}}
\def\F{{\cal F}}

\def\wt{{\rm wt\,}}

\def\deg{{\rm deg\,}}

\def\F{{\cal F}}

\def\eg{{\it e.g. }}
\def\ie{{\it i.e. }}

\def\boxit#1#2{\setbox1=\hbox{\kern#1{#2}\kern#1}%

\dimen1=\ht1 \advance\dimen1 by #1 \dimen2=\dp1 \advance\dimen2 by #1
\setbox1=\hbox{\vrule height\dimen1 depth\dimen2\box1\vrule}%
\setbox1=\vbox{\hrule\box1\hrule}%
\advance\dimen1 by .4pt \ht1=\dimen1
\advance\dimen2 by .4pt \dp1=\dimen2 \box1\relax}

\font\tensym=msbm10
\font\sevensym=msbm7
\font\fivesym=msbm5
\newfam\ssymfam
\textfont\ssymfam=\tensym
\scriptfont\ssymfam=\sevensym
\scriptscriptfont\ssymfam=\fivesym
\def\ssym{\fam\ssymfam\tensym}

\font\tengoth=eufm10
\font\sevengoth=eufm7
\font\fivegoth=eufm5

\newfam\gothfam
\textfont\gothfam=\tengoth
\scriptfont\gothfam=\sevengoth
\scriptscriptfont\gothfam=\fivegoth
\def\goth{\fam\gothfam\tengoth}

\def\mathalign#1{\hbox to 0pt{\hss$\vcenter{\openup\jot
\tabskip=0pt plus1fil \halign to \hsize
{\hfil$\displaystyle ##$\tabskip=0pt&&$\displaystyle ##$\hfil
\tabskip=0pt plus1fil \cr #1}}$\hss}}

\def\makestrut#1#2{{\dimen12=#2
\divide\dimen12 by 4\dimen11=\dimen12\multiply\dimen11 by 3
\global\setbox#1=\hbox{\vrule height\dimen11 depth\dimen12 width0pt}}}

\newdimen\tadhdimen \newdimen\tabhdimen \newdimen\vdimen
\newdimen\smtadhdimen \newdimen\smtabhdimen
\newcount\splurt
\newbox\tadstrut \newbox\tabstrut
\newbox\smtadstrut \newbox\smtabstrut

\def\setyoungsize#1#2{          
  \tadhdimen=#1\tabhdimen=#1\advance\tabhdimen by -0.4truept%
  \vdimen=#2%
  \makestrut\tadstrut\vdimen
  \advance\vdimen by -0.4pt%
  \makestrut\tabstrut\vdimen}

\def\setsmyoungsize#1#2{        
  \smtadhdimen=#1\smtabhdimen=#1\advance\smtabhdimen by -0.4truept%
  \vdimen=#2%
  \makestrut\smtadstrut\vdimen
  \advance\vdimen by -0.4pt%
  \makestrut\smtabstrut\vdimen}

\setyoungsize{15pt}{15pt}
\setsmyoungsize{9pt}{9pt}

\def\youngt#1{%
  \vcenter{\offinterlineskip
  \halign{&\copy\tadstrut\hbox to \tadhdimen{\hss$##$\hss}\cr #1}}}
\def\youngd#1{%
  \vcenter{\offinterlineskip
  \halign{&\vrule##&\copy\tabstrut\hbox to \tabhdimen{\hss$##$\hss}\cr #1}}}

\def\smyoungt#1{{\vcenter{\offinterlineskip
  \halign{&\copy\smtadstrut
  \hbox to \smtadhdimen{\hss$\scriptstyle ##$\hss}\cr #1}}}}
\def\smyoungd#1{{\vcenter{\offinterlineskip
  \halign{&\vrule##&\copy\smtabstrut
  \hbox to \smtabhdimen{\hss$\scriptstyle ##$\hss}\cr #1}}}}

\def\hdashfill{\leaders\hbox to 3pt{%
\hfil\vrule width1.5pt height0.4pt depth0pt}\hfill}

\def\smbox#1{\smyoungd{
  \multispan3\hrulefill\cr &#1&\cr\multispan3\hrulefill\cr}}

\newcommand{\bi}{\begin{itemize}}
\newcommand{\ei}{  \end{itemize}}
\newcommand{\be}{\begin{equation}}
\newcommand{\ee}{  \end{equation}}
\newcommand{\bea}{\begin{eqnarray}}
\newcommand{\eea}{  \end{eqnarray}}

\begin{document}

\title{ \bf Branching functions of $A_{n-1}^{(1)}$
       and \\
       Jantzen-Seitz 
problem for Ariki-Koike algebras }

\author{
Omar {\sc Foda}\thanks{
  Department of Mathematics, University of Melbourne,
  Parkville, Victoria 3052, Australia.}, 
Bernard {\sc Leclerc}\thanks{
  D\'epartement de Math\'ematiques, 
  Universit\'e de Caen, 
  BP 5186, 14032 Caen Cedex, France.}, 
Masato {\sc Okado}\thanks{
  Department of Mathematical Sciences, 
  Faculty of Engineering Science, 
  Osaka University, 
  Osaka 560, Japan.},\\ 
Jean-Yves {\sc Thibon}\thanks{
  Institut Gaspard Monge, 
  Universit\'e de Marne-la-Vall\'ee,
  93166 Noisy-le-Grand Cedex, France.}, 
Trevor A. {\sc Welsh}$^{\ast}$
  }

\date{October 1997}

\maketitle

\vskip 2cm
\begin{abstract}
We study the restrictions of simple modules 
of Ariki-Koike algebras $\H_m(\v)$ 
with set of parameters 
$\v= (\zeta;\zeta^{v_0},\ldots ,\zeta^{v_{l-1}})$,
where $\zeta$ is an $n$th root of unity, to their
subalgebras $\H_{m-j}(\v)$.
Using a theorem of Ariki and the crystal basis
theory of Kashiwara, we relate this problem to
the calculation of tensor product multiplicities of 
highest weight irreducible representations of 
the affine Lie algebra $A_{n-1}^{(1)}$.
These multiplicities have a combinatorial description
in terms of higher level paths or highest-lift multipartitions.

This enables us to solve the Jantzen-Seitz problem
for Ariki-Koike algebras, that is, to determine which 
irreducible representations of $\H_m(\v)$ restrict to 
irreducible representations of $\H_{m-1}(\v)$.
From a combinatorial point of view, this problem is identical
to that of computing the tensor product of 
an $A_{n-1}^{(1)}$-module of level $l$ and one of level 1.

We also consider natural generalisations of the
Jantzen-Seitz problem corresponding to the product of
a level $l$ module by a level $l'>1$ module,
and from the commutativity of tensor products,
we deduce a remarkable symmetry between 
the generalised Jantzen-Seitz conditions and the sets of parameters
of the Ariki-Koike algebras. 
 
\end{abstract}

\vfill
\eject

\setcounter{secnumdepth}{10}

\section{Introduction}

In 1994, Ariki and Koike \cite{AK} introduced an analogue of the
Iwahori-Hecke algebra for the complex reflection group
$G(l,1,m)\!=\!(\Z/l\Z)\wr \SG_m$.
This algebra $\H_m(\v)\!=\!\H_m(v;u_0,\ldots ,u_{l-1}\!)$
depends on $l+1$ parameters $v,u_0,\ldots ,u_{l-1}$
and reduces to the group algebra of $G(l,1,m)$ when
$v=1$ and $u_k=\exp(2ik\pi/l)$. 
Also, it generalises the Iwahori-Hecke algebras of the Coxeter
groups of types $A_{m-1}$ and $B_m$, which are obtained for
$l=1$ and $l=2$ respectively.
The algebra $\H_m(\v)$ appeared independently in \cite{BM}, 
where Hecke algebras of other types of complex reflection groups 
were also defined.

When the parameters are generic, $\H_m(\v)$ is semisimple
and its simple modules $S(\blambda)$
are labelled by $l$-tuples of
partitions $\blambda = (\lambda^{(0)},\ldots , \lambda^{(l-1)})$
such that $|\blambda|:=\sum_j|\lambda^{(j)}|=m$.
In this paper, we are concerned with the following choice of parameters.
Fix an integer $n\ge 2$, $l$ integers 
$0\le v_0 \le \ldots \le v_{l-1}<n$, and set
$$
v=\zeta =\exp(2i\pi/n),\qquad u_k = \zeta^{v_k}, \quad (0\le k <l).
$$
The corresponding AK-algebra will be denoted by
${\cal H}_m({\bf i})$, where ${\bf i}=(i_0,\ldots ,i_{n-1})$
and $i_k$ is the number of $v_j$ equal to $k$.
For example, ${\cal H}_m(1,0,0)$ is the Hecke algebra of type
$A_{m-1}$ at a 3rd root of 1.
It is known that ${\cal H}_m({\bf i})$ is not semisimple \cite{Ar1}.
In a recent paper \cite{Ar2} Ariki, inspired by a previous
conjecture of Lascoux, Leclerc and Thibon for type $A$ Hecke algebras
\cite{LLT1,LLT}, proved the following theorem.

Let $G_0({\cal H}_m({\bf i}))$ be the Grothendieck group of
finitely generated ${\cal H}_m({\bf i})$-mo\-dules and let 
${\cal G}({\bf i}) = \bigoplus_m G_0({\cal H}_m({\bf i}))$.
Then the infinite-dimensional complex vector space 
${\cal G}_\C({\bf i})=\C\otimes_\Z{\cal G}({\bf i})$
is naturally endowed with the structure of an $A_{n-1}^{(1)}$-module, 
isomorphic to the level $l$ irreducible representation $V(\La_{\bf i})$
of highest weight $\La_{\bf i} = \sum_k i_k\La_k$.
(Here, the $\La_k$ denote the fundamental weights of 
$A_{n-1}^{(1)}$.)
Moreover, the basis of ${\cal G}({\bf i})$ consisting of 
the classes of simple modules coincides with
Kashiwara\rq s upper global basis (or Lusztig\rq s dual canonical
basis) of $V(\La_{\bf i})$.
In particular, the simple ${\cal H}_m({\bf i})$-modules
are in one-to-one correspondence with the vertices 
of the crystal graph $B(\La_{\bf i})$ of $V(\La_{\bf i})$ that
correspond to weight vectors of principal degree $m$.
 
It turns out that this canonical basis may be efficiently computed
by embedding the $U_q(A_{n-1}^{(1)})$-module
$V_q(\La_{\bf i})$ in a suitable $q$-deformed Fock space.
In the level 1 case, this was explained in \cite{LLT}
using the Fock space construction of Misra and Miwa \cite{MM}.
The construction was generalised to arbitrary levels by Mathas \cite{Ma},
who thereby obtained the first complete classification
of simple ${\cal H}_m({\bf i})$-modules \cite{Ma2}. 
This Fock space construction yields a labelling of the vertices of
the crystal graph by either a certain class 
of $l$-tuples of Young diagrams, or by a set 
of combinatorial objects called paths, which were 
introduced in the context of solvable lattice models
by Andrews, Baxter and Forrester \cite{ABF} and extensively studied by
the Kyoto group (see \eg \cite{DJKMO}).

In this paper, we shall use Ariki's theorem to study 
the restriction of a simple ${\cal H}_m({\bf i})$-module
to the subalgebras ${\cal H}_{m-j}({\bf i})$.
In fact, we shall use some refined restriction 
operations defined by means of a natural central
element $c_m\in {\cal H}_m({\bf i})$ \cite{Ar2}.
Indeed, if $M$ is a ${\cal H}_m({\bf i})$-module,
the restricted module $M\!\downarrow_{{\cal H}_{m-1}({\bf i})}$
splits into a direct sum of eigenspaces of $c_{m-1}$,
which we denote by 
$$M\!\downarrow\!_k\,, \qquad (0\le k <n).$$
These $k$-restriction operators were first defined
for symmetric groups by Robinson (see \cite{JK}). 
More generally, we write ${M\!\downarrow \!_k}^j$
for the ${\cal H}_{m-j}({\bf i})$-module obtained from
$M$ by $j$ successive $k$-restrictions.

Let $D$ be a simple ${\cal H}_m({\bf i})$-module, and
let ${\bf j} = (j_0,\ldots ,j_{n-1})$ be another $n$-tuple of
nonnegative integers.
We say that $D$ satisfies the generalised
Jantzen-Seitz condition {\rm JS({\bf j})} if and only if
$$
{D\!\downarrow\!_k}^{j_k+1}=0 \,, \qquad (k=0,\ldots ,n-1)\,.
$$
In this case, we write $D \in {\rm JS({\bf j})}$
(or also $b\in \JS({\bf j})$ if $b$ is the vertex of $B(\La_{\bf i})$
corresponding to $D$).
This is indeed a generalisation of the original
Jantzen-Seitz condition, namely, we show
that the (ordinary) restriction
$D\!\downarrow_{{\cal H}_{m-1}({\bf i})}$ 
is irreducible if and only if 
$$
D \in \JS(1,0,\ldots ,0) \cup 
\JS(0,1,0,\ldots ,0) \cup \cdots \cup \JS(0,\ldots ,0,1)\,.
$$
Let $b$ be the vertex of the crystal graph of $V(\La_{\bf i})$
corresponding to $D$.
It follows from Ariki's theorem that $D\in \JS({\bf j})$
if and only if 
$$
\tilde e_k^{j_k+1}\, b = 0 \,, \qquad (k=0,\ldots ,n-1)\,,
$$
where $\tilde e_k$ denotes the Kashiwara operator acting 
on the crystal graph.
This allows us to relate the generalised JS-problem to
the calculation of tensor product multiplicities.
Write $\La_{\bf j} = \sum_k j_k\La_k$, 
and denote by $\deg b$ the homogeneous degree of $b\in B(\La_{\bf i})$. 
Then we obtain
$$
\sum_{b\in B(\La_{\bf i})\cap \JS({\bf j})}
z^{\deg b}
=
\sum_\La b^{\Lambda}_{\Lambda_{\bf j},\Lambda_{\bf i}}(z)\,,
$$ 
where the 
$ b^{\Lambda}_{\Lambda_{\bf j},\Lambda_{\bf i}}(z)$
are the branching functions of the tensor
product $V(\La_{\bf j})\otimes V(\La_{\bf i})$.
Thus, the generating function of the number of simple $\H_m({\bf i})$-modules
satisfying the condition $\JS({\bf j})$ is a sum of branching functions
of $A_{n-1}^{(1)}$.

One striking fact about this formula is that ${\bf i}$
and ${\bf j}$ play symmetric roles, which means that
the number of 
simple ${\cal H}_m({\bf i})$-modules satisfying the
condition $\JS({\bf j})$ is the  same as the 
number of 
simple ${\cal H}_m({\bf j})$-modules satisfying the 
condition $\JS({\bf i})$.

This result motivates a deeper study of the
sets of multipartitions 
${\cal Y}(\Lambda_{\bf j},\Lambda_{\bf i})$
which label the highest weight vectors of 
$V(\La_{\bf j})\otimes V(\La_{\bf i})$.
When both $\Lambda_{\bf j}$ and $\Lambda_{\bf i}$ are
of level 1, \ie when $\Lambda_{\bf j} = \Lambda_j$
and $\Lambda_{\bf i} = \Lambda_i$ are fundamental weights,
the set ${\cal Y}(\Lambda_j,\Lambda_i)$
consists of a class of partitions described by
Foda, Okado and Warnaar \cite{FOW} in their
work on RSOS solvable lattice models based
on the coset $(\slchap_n)_1 \times (\slchap_n)_1 /
(\slchap_n)_2$.
The surprising observation that these partitions are the same as
those occurring in the paper of Jantzen and Seitz \cite{JS}
in the modular representation theory of symmetric groups 
was in fact the starting point of this work \cite{FLOTW1}.

In this article we describe, much in the spirit
of \cite{FOW}, the elements of
${\cal Y}(\La_{\bf j},\La_{\bf i})$ in the case
when $\La_{\bf j}$ is of level one and $\La_{\bf i}$ is arbitrary.
To this aim, we first consider the
problem of relating the two sets ${\cal Y}(\La_{\bf j},\La_{\bf i})$
and ${\cal Y}(\La_{\bf i},\La_{\bf j})$ for arbitrary dominant
integral weights $\La_{\bf i}$ and $\La_{\bf j}$.
Using the combinatorics of paths, we construct a simple bijection
between ${\cal Y}(\La_{\bf j},\La_{\bf i})$ and 
${\cal Y}(\sharp\La_{\bf i},\sharp\La_{\bf j})$, where
$\sharp$ denotes the root diagram automorphism exchanging 
$\La_i$ and $\La_{n-i}$.
We then apply this construction in the case of a fundamental weight
$\La_{\bf j} = \La_j$ and thereby obtain the desired characterisation.

The paper is divided into two main parts.
The first one, Section~2, is devoted to the calculation of branching functions
using the combinatorics of paths and multipartitions.
Since this language is probably not familiar to many readers
interested in AK-algebras, we feel it appropriate 
to explain it in some detail.
Thus, Sections~2.1 and 2.2 are almost entirely expository, our
main references being \cite{DJKMO} and \cite{JMMO}.
In 2.1, we review the Fock space representations and their
$q$-deformations, and describe how
the vertices of the crystal graph of $V(\La)$ may be labelled using
either unrestricted paths or their corresponding highest-lift multipartitions.
In 2.2, we use crystal basis theory to identify some specific
classes of paths, called restricted paths, whose generating
functions are the branching functions.
We also introduce the corresponding restricted multipartitions.
Section~2.3 deals with the involution $\sharp$ and 
Section~2.4 contains the description of ${\cal Y}(\La_j,\La_{\bf i})$
obtained via $\sharp$.

The second part, Section~3, deals with AK-algebras.
Section~3.1 introduces, following \cite{AK,Ar1,Ar2},
the Ariki-Koike algebra and 
provides a realisation of it as a quotient of the affine 
Hecke algebra.
Section~3.2 recalls the basic tools and results of the representation
theory of AK-algebras, including
the definition of the $k$-restriction operators and the
statement of Ariki's theorem.
Finally, in Sections~3.3 and 3.4, we present our results on the
generalised JS-problem and interpret
the combinatorial theorems of 2.3 and 2.4 in this context.

\section{Combinatorics of branching functions of $A_{n-1}^{(1)}$}

\subsection{Highest weight modules, unrestricted paths and highest-lifts}
\label{SECT2.1}

We begin by recalling some standard notation (see {\em e.g.} 
\cite{KacBook} \cite{JMMO}).

Let $e_i,f_i,h_i, (0\le i \le n-1)$ be the Chevalley
generators of the affine Lie algebra $\slchap_n=A_{n-1}^{(1)}$.
The degree generator and the canonical central element 
are denoted respectively by $d$ and $c$.
The derived algebra $\slchap_n' = [\slchap_n , \slchap_n]$
is the subalgebra obtained by omitting the degree generator $d$.

Let $\Lambda_0, \ldots , \Lambda_{n-1}$ denote the
fundamental weights and $\delta$ the null root.
For $i\in \Z$, set $\Lambda_{i}=\Lambda_{(i\mod n)}$,
$\alpha_i=2\Lambda_i-\Lambda_{i-1}-\Lambda_{i+1}
+ \delta_{i0}^{(n)} \,\delta$ 
(where $\delta_{ij}^{(n)}=1$ if $(i-j)\mod n=0$ 
and $\delta_{ij}^{(n)}=0$ otherwise)
and $\epsilon_i=\Lambda_{i+1}-\Lambda_i$.
The $\alpha_i$ are the simple roots and
the $\epsilon_i$ can be identified with the weights
of the $n$-dimensional defining representation of
$\Sl_n$ extended to an $\slchap_n'$-representation.
%
The root lattice and the weight lattice are respectively 
$Q=\bigoplus_{i=0}^{n-1} \Z\alpha_i$
and
$P=\Z\delta \oplus (\bigoplus_{i=0}^{n-1} \Z\Lambda_i)$. 
For $l\in\N$, let 
$$
P_l=\left\{\sum_{i=0}^{n-1}a_i\Lambda_i\mid a_i\in\Z,
\ \sum_{i=0}^{n-1}a_i=l\right\},\qquad
P_l^+=\left\{\sum_{i=0}^{n-1}a_i\Lambda_i\in P_l\mid a_i\ge0\right\}
$$
be the set of (classical) level $l$ integral weights and
the subset of dominant ones, respectively.
Finally, let 
$
P^+=\Z\delta \oplus( \bigoplus_{i=0}^{n-1} \N\Lambda_i)
$
be the set of dominant integral weights.

For each $\Lambda \in P^+$ there exists a unique integrable 
highest weight module $V(\Lambda)$ with highest weight $\Lambda$.
It can be explicitly constructed as a submodule of a Fock space
$\F(\Lambda)$ \cite{DJKMO}, as we now explain.

Throughout this section we fix $\Lambda \in P_l^+$ and 
we write
\[
\Lambda = \Lambda_{v_0}+\Lambda_{v_1}+\cdots+\Lambda_{v_{l-1}}\,,
\]
where we may assume that 
$0\le v_0\le v_1\le\cdots\le v_{l-1}<n$.
The level $l$ Fock space $\F(\Lambda)$ is a tensor 
product of $l$ level 1 Fock spaces
$$
\F(\Lambda) = \F(\Lambda_{v_0}) \otimes \cdots \otimes
\F(\Lambda_{v_{l-1}})\,,
$$
and each level 1 Fock space representation is constructed as follows.
As a vector space $\F(\Lambda_j) = \bigoplus_\lambda \C\,  u_\lambda$,
where $\lambda$ runs through the set $\Pi$ of all partitions.
In other words, $\F(\Lambda_j)$ is an infinite dimensional
$\C$-space with a distinguished basis $u_\lambda$ labelled
by Young diagrams.
To define the action of $\slchap_n$ on $\F(\Lambda_j)$, one 
uses a colouring of the nodes of each $\lambda \in \Pi$,
namely the node in the $i$th row and the $k$th column of
$\lambda$ is filled with the colour $(k-i+j)\mod n$.
We write $\mu/\lambda=\smbox{r}$
to indicate that the Young diagram $\mu$ is obtained 
from $\lambda$ by adding a node with colour $r$.
Then $\slchap_n$ acts on $\F(\Lambda_j)$ by \cite{DJKMO}
\begin{eqnarray}
f_r u_\lambda
&=& \sum_{\mu/\lambda=\smbox{r}} u_\mu\,,  \label{EQ1}\\
e_r u_\mu
&=& \sum_{\mu/\lambda=\smbox{r}} u_\lambda\,, \label{EQ2}\\
c \,u_\lambda & = & u_\lambda\,, \label{EQ3}\\
d \,u_\lambda &=& -N^0(\lambda) u_\lambda \,,\label{EQ4}
\end{eqnarray} 
where $N^0(\lambda)$ is the number of $0$-nodes of $\lambda$.   

The level $l$ action of $\slchap_n$ on $\F(\Lambda)$ is  
obtained by taking the tensor product of these level
1 actions in the standard way.
The natural basis of $\F(\Lambda)$ consists of the monomial
tensors
$$
u_\blambda = u_{\lambda^{(0)}}\otimes \cdots \otimes u_{\lambda^{(l-1)}}
$$
indexed by $l$-tuples of partitions 
$\blambda=(\lambda^{(0)},\ldots,\lambda^{(l-1)})$.
 
\begin{example} \label{EX1} {\rm
Take $n=3$ and $\Lambda = 2\Lambda_1 + \Lambda_2$.
Set $\blambda = ((3,2),(1,1,1),(5,4,1))$. 
The coloured Young diagram of $\blambda$ is
$$
\smyoungd{\omit\phantom{\hrulefill}\cr
 \omit&\cr
 \multispan{8}\hrulefill\cr
 &1&&2&&0&&\cr
 \multispan{7}\hrulefill\cr
 &0&&1&&\cr
 \multispan{5}\hrulefill\cr
 &\cr
 \multispan{1}\hrulefill\cr
 &\cr
 \multispan{1}\hrulefill\cr}
\quad
\smyoungd{\omit\phantom{\hrulefill}\cr
 \omit&\cr
 \multispan{4}\hrulefill\cr
 &1&&\cr
 \multispan{3}\hrulefill\cr
 &0&&\cr
 \multispan{3}\hrulefill\cr
 &2&&\cr
 \multispan{3}\hrulefill\cr
 &\cr
 \multispan{1}\hrulefill\cr}
\quad
\smyoungd{\multispan{12}\hrulefill\cr
 &2&&0&&1&&2&&0&&\cr
 \multispan{11}\hrulefill\cr
 &1&&2&&0&&1&&\cr
 \multispan{9}\hrulefill\cr
 &0&&\cr
 \multispan{3}\hrulefill\cr
 &\cr
 \multispan{1}\hrulefill\cr
 &\cr
 \multispan{1}\hrulefill\cr}
$$
The lowering operators of $\slchap_3$ act on 
$u_\blambda \in {\cal F}(\Lambda)$ as follows:
\begin{eqnarray*}
f_0 u_\blambda & = & 0, \\ 
f_1 u_\blambda & = & u_{(({\bf 4},2),(1,1,1),(5,4,1))}
+  u_{((3,2),(1,1,1,{\bf 1}),(5,4,1))}
+  u_{((3,2),(1,1,1),({\bf 6},4,1))} \\
&& +\  u_{((3,2),(1,1,1),(5,4,{\bf 2}))}, \\
f_2 u_\blambda & = & u_{((3,{\bf 3}),(1,1,1),(5,4,1))}
+ u_{((3,2,{\bf 1}),(1,1,1),(5,4,1))}
+  u_{((3,2),({\bf 2},1,1),(5,4,1))} \\
&&+\  u_{((3,2),(1,1,1),(5,{\bf 5},1))} 
+  u_{((3,2),(1,1,1),(5,4,1,{\bf 1}))}. 
\end{eqnarray*}
(The parts of $\blambda$ which have been increased by 1
are printed in bold type).
The degree operator acts by $d\, u_\blambda = -7 u_{\blambda}$.
\endex}
\end{example}
%
It follows from (\ref{EQ1}), (\ref{EQ2}), (\ref{EQ3}) and (\ref{EQ4}),
that for all $\blambda$, $u_\blambda$ is a weight vector
of $\F(\Lambda)$ of weight $\wt(\blambda)$, where we define
\begin{equation}\label{MultiWtDef}
\wt(\blambda) = \Lambda - \sum_{i=0}^{n-1} N^i(\blambda) \alpha_i\,,
\end{equation}
and where
$N^r(\blambda)$ is the total number of $r$-nodes of $\blambda$.

Let $u_\emptyset$ denote the vacuum vector of $\F(\Lambda)$,
that is, the vector labelled by the empty multipartition.
It is clear that $u_\emptyset$
is a highest weight vector of $\F(\Lambda)$ of weight $\Lambda$.
The submodule $U(\slchap_n)\,u_\emptyset$ is isomorphic to $V(\Lambda)$.
In  \cite{DJKMO}, an explicit basis of this submodule was constructed
using combinatorial objects called {\em paths}.
These we now describe.

Let 
$
{\cal A}_l^+ = \{ \sum_{i=0}^{n-1} a_i \epsilon_i \ |\
a_i \in \N, \sum a_i = l \}
$
denote the set of weights of the $l$th symmetric power of $\C^n$
(considered as an $\slchap_n'$-representation).
A {\em level $l$ path} is an infinite sequence $p=(p_0,p_1,p_2,\ldots)$
of weights $p_k\in P_l$.
A {\em $\Lambda$-path} is a level $l$ path such that:
\begin{eqnarray}
&& \mbox{for all $k\ge0$, $p_{k+1}-p_k \in {\cal A}_l^+$,}\\ 
&& \mbox{for $k$ large enough, $p_k =
\Lambda_{v_0+k}+\cdots+\Lambda_{v_{l-1}+k}$.}\label{GSP}
\end{eqnarray}
Condition~(\ref{GSP}) means that, 
except for a finite number of indices, 
$p= (p_k)$ coincides with the {\em ground state path} 
$$\overline p = (\overline p_k) =
(\Lambda_{v_0+k}+\cdots+\Lambda_{v_{l-1}+k})\,.$$
The smallest integer $k_*$ such that $p_k = \overline p_k$
for $k\ge k_*$ is called the {\em length} of $p$ and is 
denoted by $\ell(p)$.
The set of $\Lambda$-paths is denoted by ${\cal P}(\Lambda)$.

Note that a $\Lambda$-path $p=(p_k)$ is completely determined
by the associated sequence $\eta = (\eta_k)$ of elements 
$\eta_k = p_{k+1} - p_k$
of ${\cal A}_l^+$. 
It is sometimes more convenient to think of a $\Lambda$-path
in this way as a sequence $\eta$ of elements
of ${\cal A}_l^+$ with the tail condition 
$$
\eta_k = \overline \eta_k
:= \epsilon_{v_0+k}+\cdots+\epsilon_{v_{l-1}+k} \,,
\quad (k>\!>0) \,.
$$
We shall sometimes write
$\eta_k = \sum_{j=0}^{l-1}\epsilon_{\gamma_j(k)}\,$,
where $0\le\gamma_j(k)<n$ for $k=0,1,\ldots,l-1$.

The notion of path comes from the study of solvable lattice 
models, where, roughly speaking, a path corresponds to an eigenvector
of the corner transfer matrix at the absolute temperature $q=0$
(see {\em e.g.} \cite{FLOTW2}).
It was found by Date, Jimbo, Kuniba, Miwa and Okado \cite{DJKMO} that
the set of $\Lambda$-paths is an appropriate labelling
for a basis of weight vectors of $V(\Lambda)$.
Namely, in the context of statistical mechanics, a path $p$
has an {\em energy} $E(p)$ defined as follows.
One introduces a function $H$ on ${\cal A}_l^+ \times {\cal A}_l^+$
given, for $\alpha = \epsilon_{\mu_0}+\cdots +\epsilon_{\mu_{l-1}}$
and $\beta = \epsilon_{\nu_0}+\cdots +\epsilon_{\nu_{l-1}}$,
by 
\begin{equation}\label{HDEF}
H(\alpha,\beta) = \min_{\sigma \in \SG_l}\, 
\sum_{i=0}^{l-1} \theta(\mu_i - \nu_{\sigma(i)}) \,,
\end{equation}
where $\theta(a)=1$ if $a\ge 0$ and $\theta(a)=0$ otherwise.
Thus, for $n=3=l$, one has
$$
H(\epsilon_0+\epsilon_1+\epsilon_2,
\epsilon_0+\epsilon_1+\epsilon_2) = \theta(0-1) + \theta(1-2) +
\theta(2-0) = 1\,.
$$
Then one puts
\begin{equation}\label{EDEF}
E(p) = \sum_{k=1}^\infty
k\left( H(\eta_{k-1},\eta_k) - H(\overline \eta_{k-1}, \overline
\eta_k)\right)
\,.
\end{equation}
This is in fact a finite sum
since, by condition~(\ref{GSP}) above, 
$\eta_k = \overline \eta_k$ for $k$ large enough.
Finally one defines the {\em weight} $\wt(p)\in P$ of $p$ by:
\begin{equation}\label{PathWtDef}
\wt(p) = p_0 - E(p)\delta\,.
\end{equation}

\begin{theorem}{\rm \cite{DJKMO}}\label{TH2.2}
Let $\Lambda \in P_l^+$.
The formal character of the representation $V(\Lambda)$ of $\slchap_n$
is given by
$$
{\rm ch} V(\Lambda) = \sum_{p\in {\cal P}(\Lambda)} e^{\wt(p)} \,.
$$
\end{theorem}
\begin{figure}[t]
\begin{center}
\leavevmode
\epsfxsize =13.5cm
\epsffile{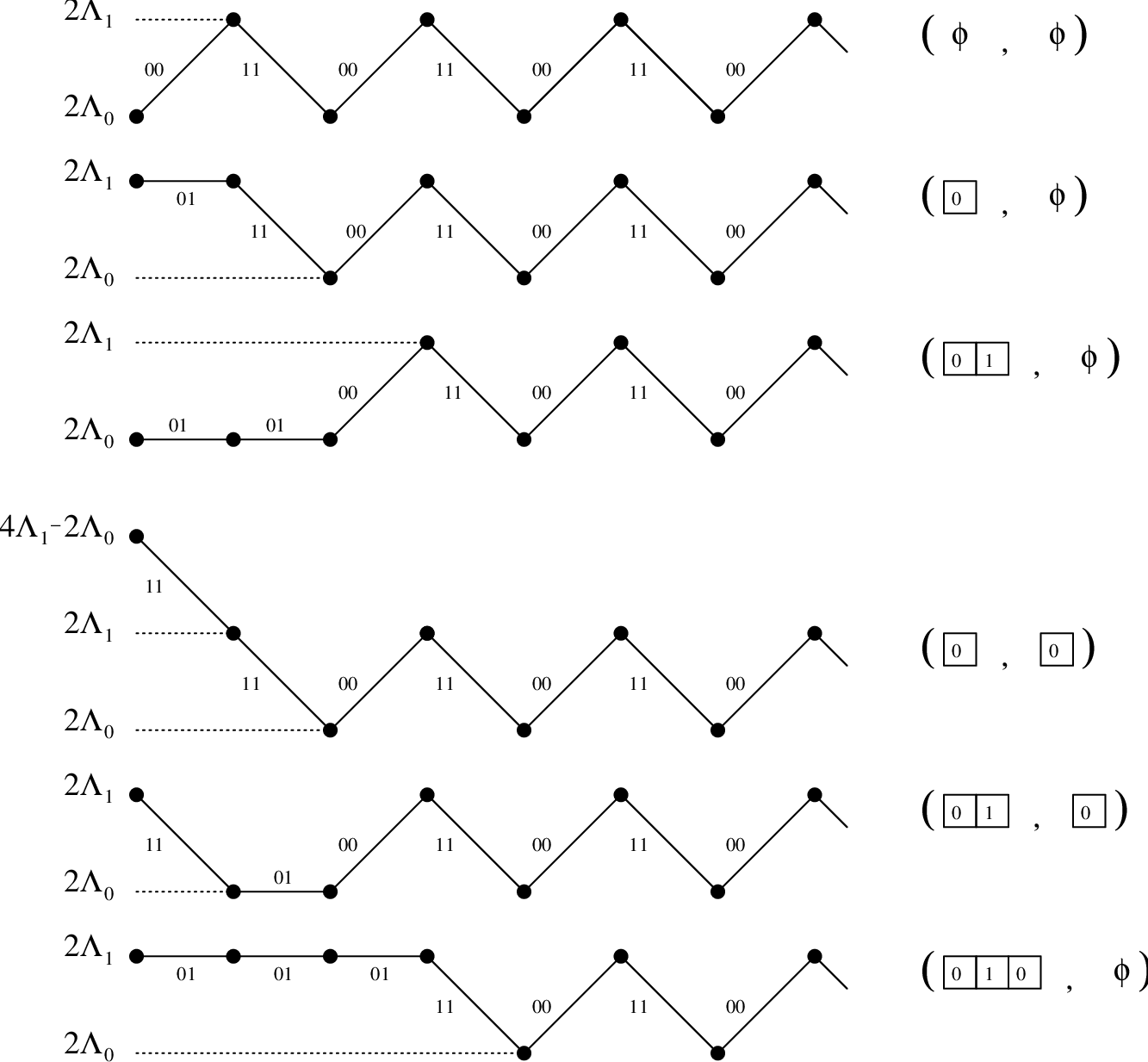}
\end{center}
\caption{\label{FIG0} $2\La_0$-paths for $\slchap_2$ 
and their highest-lift bipartitions}
\end{figure}
\begin{example}{\rm
We illustrate the theorem by enumerating a few paths in 
${\cal P}(2\Lambda_0)$ for $\slchap_2$.
We write for short $00$, $01$, $11$, in place of
$2\epsilon_0$, $\epsilon_0 + \epsilon_1$, $2\epsilon_1$.
The paths are given as sequences $(\eta_k)$ of elements
of ${\cal A}_2^+$.
$$
\matrix{
\hbox{\rm path}       & \hbox{\rm energy}  & \hbox{weight} \cr
(00,11,00,11,\ldots ) & 0                  & 2\Lambda_0    \cr
(01,11,00,11,\ldots ) & 1                  & 2\Lambda_1 - \delta   \cr
(01,01,00,11,\ldots ) & 1                  & 2\Lambda_0 - \delta   \cr
(11,11,00,11,\ldots ) & 2                  & 4\La_1-2\La_0 - 2\delta  \cr
(11,01,00,11,\ldots ) & 2                  & 2\Lambda_1 - 2\delta  \cr
(01,01,01,11,\ldots ) & 2                  & 2\Lambda_1 - 2\delta  
}
$$
These paths are shown in Figure~\ref{FIG0}, with the
corresponding highest-lift multipartitions which will
be defined below.
Thus, one has
$$
{\rm ch} V(2\Lambda_0) = e^{2\Lambda_0} + e^{2\Lambda_1 - \delta }
+ e^{2\Lambda_0 - \delta } + e^{ 4\La_1-2\La_0 - 2\delta }
+2e^{2\Lambda_1 - 2\delta} + \cdots
$$
\endex}
\end{example}
The relationship between paths and highest weight representations
of $\slchap_n$ was later clarified using the crystal
basis theory of Kashiwara \cite{MM,JMMO,KMN1,KMN2}.
This involves a $q$-deformation of $\slchap_n$.

Let $U_q(\slchap_n)$ be the quantized enveloping algebra
of $\slchap_n$. We denote by $V_q(\Lambda)$ the irreducible
$U_q(\slchap_n)$-module with highest weight $\Lambda$.
We shall follow \cite{JMMO} and construct $V_q(\Lambda)$
as a submodule of a $q$-deformed level $l$ Fock space 
${\cal F}_q(\Lambda)$.
As a vector space,
$$
{\cal F}_q(\Lambda) = \sum_{\blambda \in \Pi^l}
\Q(q) v_{\blambda} \,,
$$
that is, ${\cal F}_q(\Lambda)$ has a distinguished $\Q(q)$-basis
$\{v_\blambda\}$ labelled by the set of all $l$-tuples of
partitions 
$\blambda=(\lambda^{(0)},\ldots,\lambda^{(l-1)})$.

To describe the action of $U_q(\slchap_n)$ requires some
notation.
First we colour the nodes of $\blambda$ as before, filling
the node of $\lambda^{(j)}$ which lies on the $i$th row
and the $k$th column with $r=(k-i+v_j) \mod n$.
We write $\bmu/\blambda=\smbox{r}$
to indicate that the multipartition $\bmu$ is obtained
from $\blambda$ by adding a node with colour $r$.
In this case, we say that $\bmu/\blambda$ is a {\em removable $r$-node}
of $\bmu$ or an {\em addable $r$-node} of $\blambda$.
With each removable or addable $r$-node of $\blambda$,
there is associated a pair $(d,j)$ of integers, where 
$d=k-i+v_j$ indicates the actual diagonal containing that node,
and $j$ indicates the component $\lambda^{(j)}$ of 
$\blambda$ to which the node belongs.
Then, define a total order on the set of removable and addable
$r$-nodes of $\blambda$ by:
\begin{equation}
(d,j) < (d',j') \qquad \Longleftrightarrow \qquad
\left( (d<d') \hbox{ \rm or } (d=d' \hbox{ \rm and } j > j')\right)\,.
\end{equation}
\begin{example}\rm\label{EX2}
Take $n=3$, $l=2$, $\Lambda = \Lambda_0 + \Lambda_1$
and consider the bipartition 
$$\blambda = ((9,8,7,5,4,4,1,1),(9,9,7,6,5,3,3))$$
with coloured diagram
$$
\smyoungd{\omit\phantom{\hrulefill}\cr
 \omit&\cr
 \multispan{20}\hrulefill\cr
 &0&&1&&2&&0&&1&&2&&0&&1&&2&&\cr
 \multispan{19}\hrulefill\cr
 &2&&0&&1&&2&&0&&1&&2&&0&\cr
 \multispan{17}\hrulefill\cr
 &1&&2&&0&&1&&2&&0&&1&\cr
 \multispan{15}\hrulefill\cr
 &0&&1&&2&&0&&1&\cr
 \multispan{11}\hrulefill\cr
 &2&&0&&1&&2&\cr
 \multispan{9}\hrulefill\cr
 &1&&2&&0&&1&\cr
 \multispan{9}\hrulefill\cr
 &0&\cr
 \multispan{3}\hrulefill\cr
 &2&\cr
 \multispan{3}\hrulefill\cr
 &\cr
 \multispan{1}\hrulefill\cr
 }
\quad
\smyoungd{\multispan{20}\hrulefill\cr
 &1&&2&&0&&1&&2&&0&&1&&2&&0&&\cr
 \multispan{19}\hrulefill\cr
 &0&&1&&2&&0&&1&&2&&0&&1&&2&&\cr
 \multispan{19}\hrulefill\cr
 &2&&0&&1&&2&&0&&1&&2&\cr
 \multispan{15}\hrulefill\cr
 &1&&2&&0&&1&&2&&0&\cr
 \multispan{13}\hrulefill\cr
 &0&&1&&2&&0&&1&\cr
 \multispan{11}\hrulefill\cr
 &2&&0&&1&\cr
 \multispan{7}\hrulefill\cr
 &1&&2&&0&\cr
 \multispan{7}\hrulefill\cr
 &\cr
 \multispan{1}\hrulefill\cr
 &\cr
 \multispan{1}\hrulefill\cr
 &\cr
 \multispan{1}\hrulefill\cr
 }
$$
The addable and removable 1-nodes of $\blambda$ are ordered in 
the following way:
$$
A_{-8,0}< A_{-5,0}< R_{-2,0}< R_{1,1}< R_{1,0}< A_{4,1}
< R_{4,0}< A_{7,0} < A_{10,1}\,,
$$
where $A_{d,j}$ means an addable node in position $(d,j)$
and $R_{d,j}$ a removable node in position $(d,j)$.
\endex
\end{example}
Now suppose that $\bmu/\blambda$ is the $r$-node $(d,j)$.
We put 
\begin{verse}
$N_r^{<}(\blambda,\bmu) = \sharp \{$ addable $r$-nodes $(d',j')$
of $\blambda$ such that $(d',j')<(d,j)\}$\\
$\qquad\quad\quad - \ \sharp \{$ removable $r$-nodes $(d',j')$ of $\blambda$
such that $(d',j')<(d,j)\}$,

$N_r^{>}(\blambda,\bmu) = \sharp \{$ addable $r$-nodes $(d',j')$
of $\blambda$ such that $(d',j')>(d,j)\}$\\
$\qquad\quad\quad - \ \sharp \{$ removable $r$-nodes $(d',j')$ of $\blambda$
such that $(d',j')>(d,j)\}$,

$N_r(\blambda) = \sharp \{$ addable $r$-nodes of $\blambda\ \}
- \ \sharp \{$ removable $r$-nodes of $\blambda\ \}$.

\end{verse}
\begin{theorem}{\rm\cite{JMMO}}\label{TH2-5}
$U_q(\slchap_n)$ acts on ${\cal F}_q(\Lambda)$ by

\smallskip
$f_r v_\blambda = \sum_\bmu q^{N_r^{>}(\blambda,\bmu)} \, v_\bmu \,,$
sum over all $\bmu$ such that $\bmu/\blambda=\smbox{r}$\,,

$e_r v_\bmu = \sum_\blambda q^{-N_r^{<}(\blambda,\bmu)} \, v_\blambda \,,$ 
sum over all $\blambda$ such that 
$\bmu/\blambda=\smbox{r}$\,,
 
\smallskip
$q^{h_r} \, v_\blambda = q^{N_r(\blambda)} \, v_\blambda\,,$

\smallskip
$q^c \,v_\blambda = q^l\, v_\blambda \,,$

\smallskip
$q^d \, v_\blambda = q^{-N^0(\blambda)} \, v_\blambda\,.$
\end{theorem}
\begin{example} \label{EX3} {\rm
We give the $q$-deformation of the calculation of
Example~\ref{EX1}.
Take $n=3$ and $\Lambda = 2\Lambda_1 + \Lambda_2$.
Set $\blambda = ((3,2),(1,1,1),(5,4,1))$. Then the
lowering operators of $U_q(\slchap_3)$ act on
$v_\blambda \in {\cal F}_q(\Lambda)$ as follows:
\begin{eqnarray*}
f_0\, v_\blambda & = & 0 \\
f_1\, v_\blambda & = & q \,v_{((4,2),(1,1,1),(5,4,1))}
+ q\, v_{((3,2),(1,1,1,1),(5,4,1))} \\
&& +\  v_{((3,2),(1,1,1),(6,4,1))} 
+   v_{((3,2),(1,1,1),(5,4,2))} \\
f_2\, v_\blambda & = & q\, v_{((3,3),(1,1,1),(5,4,1))}
+ q^3\, v_{((3,2,1),(1,1,1),(5,4,1))} \\
&&+\ q^2\, v_{((3,2),(2,1,1),(5,4,1))} 
+   v_{((3,2),(1,1,1),(5,5,1))} \\
&& +\ q^3\,  v_{((3,2),(1,1,1),(5,4,1,1))}
\end{eqnarray*}
\endex}
\end{example}
\begin{note} \label{NoteFockSpace}{\rm
It is possible, and perhaps more natural, to define
another level $l$ $q$-deformed Fock space 
${\cal F}'_q(\Lambda)$ by imitating the two steps
of the $q=1$ construction.
First, one constructs the level 1 $q$-Fock spaces
${\cal F}_q(\Lambda_{v_j})$ as above. Then, one
defines via the comultiplication of $U_q(\slchap_n)$
$$
\F_q'(\Lambda) = \F_q(\Lambda_{v_0}) \otimes \cdots \otimes
\F_q(\Lambda_{v_{l-1}})\,.
$$
Clearly, the two Fock spaces $\F_q(\Lambda)$ and $\F_q'(\Lambda)$
coincide at $q=1$. However, they are different in general.
More precisely, $\F_q(\Lambda)$ and $\F_q'(\Lambda)$
are isomorphic as $U_q(\slchap_n)$-modules but the
formulae giving the action of the Chevalley generators on
multipartitions are not the same.
This leads to two different ways of labelling the vertices
of the same crystal graph. 
In this paper, we work exclusively with $\F_q(\Lambda)$
in contrast with  Mathas \cite{Ma} who
uses $\F_q'(\Lambda)$.
As a consequence, our labelling of the simple modules of
the Ariki-Koike algebras differs from that
of \cite{Ma} (see below, Section~\ref{SECT3.4}). 
}
\end{note}

Using this explicit description of the $q$-deformed Fock space, 
it is possible to compute its crystal basis \cite{JMMO}.
To explain how the arrows of the crystal graph are obtained,
we introduce the notion of a {\em good node}.
Fix a residue $r$ between $0$ and $n-1$ and consider the
sequence of removable and addable $r$-nodes of $\blambda$,
ordered as explained above. 
This sequence is called the {\em $r$-signature} of $\blambda$.
Thus for $n=3$, $\Lambda = \Lambda_0+\Lambda_1$ and
the bipartition $\blambda$ of Example~\ref{EX2},
we get the following 1-signature:
$$
A_{-8,0} A_{-5,0} R_{-2,0} R_{1,1} R_{1,0} A_{4,1}
 R_{4,0} A_{7,0}  A_{10,1}\,.
$$
Now recursively remove $RA$ pairs (together with their
subscripts) from this sequence until none remain.
The sequence will then be of the form 
$$A\cdots AR\cdots R.$$
The remaining nodes are called the {\em normal $r$-nodes}
of $\blambda$.
The node corresponding to the leftmost $R$ is termed a
{\em good removable $r$-node} and that corresponding to the
rightmost $A$ is termed a {\em good addable $r$-node} (note that
there is at most one of each).
In the example above, the removal procedure first disposes of $R_{1,0}
A_{4,1}$ and
$R_{4,0} A_{7,0}$, and then disposes of $R_{1,1} A_{10,2}$ so that
$$
A_{-8,0} A_{-5,0} R_{-2,0}
$$
remains. Therefore $\blambda$ has a good addable node on the $-5$
diagonal of $\lambda^{(0)}$
and a good removable node on the $-2$ diagonal of $\lambda^{(1)}$.
In the case of 0-nodes and 2-nodes, we first obtain the sequences
\begin{eqnarray*}
&&A_{-6,1} R_{-3,1} A_{0,0} R_{3,1} A_{6,1} R_{6,0} A_{9,0}\\
\noalign{\noindent\hbox{and}}
&&R_{-7,0} A_{-1,1} A_{2,1} A_{2,0} R_{5,1} A_{5,0} R_{8,1} R_{8,0}
\end{eqnarray*}
respectively which, after the removal procedure, produce
\begin{eqnarray*}
&&A_{-6,1}\\
\noalign{\noindent\hbox{and}}
&&A_{2,1} A_{2,0} R_{8,1} R_{8,0}
\end{eqnarray*}
respectively.
Thus $\blambda$ has a good addable 0-node on the $-6$ diagonal of
$\lambda^{(1)}$, but no good removable 0-node,
and a good addable 2-node on the $2$ diagonal of $\lambda^{(0)}$ 
and a good removable 2-node on the $8$ diagonal of $\lambda^{(1)}$. 

\begin{theorem}{\rm\cite{JMMO}}
Let $A\subset \Q(q)$ be the ring of rational functions without
pole at $q=0$.
Let $L = \bigoplus_{\blambda \in \Pi^l}
A\,v_\blambda$, and let $B$ be the $\Q$-basis of $L/qL$ given
by $B=\{v_\blambda \mod qL \ | \ \blambda \in \Pi^l \}$.
Then $(L,B)$ is a crystal basis of
${\cal F}_q(\Lambda)$.
Moreover, the crystal graph contains the arrow
$$
v_\blambda \stackrel{r}{\longrightarrow} v_\bmu
$$
if and only if $\bmu$ is obtained from $\blambda$ by adding
a good $r$-node.
\end{theorem}

Let $v_{\emptyset}$ denote the vacuum vector of ${\cal F}_q(\Lambda)$.
The submodule $U_q(\slchap_n) v_{\emptyset}$ is known
to be isomorphic to $V_q(\Lambda)$.
Therefore, using properties of crystal bases one can
obtain the crystal graph of $V_q(\Lambda)$ as the connected
component of that of ${\cal F}_q(\Lambda)$ which contains
$v_{\emptyset}$.
The multipartitions labelling the vertices of this subgraph
have been described in \cite{JMMO} using the following map
from multipartitions to paths.

Let $\blambda=(\lambda^{(0)},\ldots,\lambda^{(l-1)})$ and 
let $k_* = \max \{ \lambda^{(0)}_1,\ldots,\lambda^{(l-1)}_1\}$
be the greatest part of the partitions $\lambda^{(j)}$.
We define the path $p=\pi(\blambda) \in {\cal P}(\Lambda)$
of length $\ell(p)=k_*$ by
\begin{equation}\label{ETADEF}
\eta_k = \sum_{j=0}^{l-1} \epsilon_{v_j+k-{\lambda^{(j)}}_{k+1}'}\,,
\end{equation}
where ${\lambda^{(j)}}_k'$ denotes the height of the 
$k$th column of the Young diagram $\lambda^{(j)}$.
The sequence $(\eta_k)$ is easily read from the coloured 
diagram of $\blambda$ by taking
the colours of the nodes immediately under the bottom end
of each column. 
For instance, if $\blambda$ is again the bipartition $\blambda$
of Example~\ref{EX2} then, for convenience, writing
$ij$ instead of $\epsilon_i+\epsilon_j$, we have
$$
\eta = (01,11,22,02,00,22,01,02,11,01,12,02,01,12,\ldots )\,.
$$
It was proved in \cite{DJKMO,JMMO} that the map $\pi$
is surjective and that each path $p$ has a distinguished
preimage $\blambda$ called the {\em highest-lift} of $p$.
It is characterised as follows. 

We say that a multipartition $\blambda$ is {\em cylindrical
of highest weight $\Lambda$} if 
\begin{eqnarray*}
&&\lambda^{(j)}_i\ge\lambda^{(j+1)}_{i+v_{j+1}-v_j},\qquad0\le j\le l-2,
\quad i=1,2,\ldots ,\\
&&\lambda^{(l-1)}_i\ge\lambda^{(0)}_{i+n+v_0-v_{l-1}},\qquad
i=1,2,\ldots.
\end{eqnarray*}
If we make the convention to align the partitions
so that the first rows beginning with a node of colour 0 are
adjacent to one another, then the first condition
means that the lengths of adjacent rows should not
increase from left to right across the diagram. 
The second condition
can be similarly checked by putting a copy of $\lambda^{(0)}$
on the right, but raised $n$ rows.

Given two multipartitions $\blambda$, $\bmu$, we say that
$\blambda$ is {\em higher than} $\bmu$ if
$$
{\lambda^{(j)}}'_k\le{\mu^{(j)}}'_k,\qquad0\le j\le l-1,
\quad k=1,2,\ldots ,
$$
that is, if the frontier of the Young diagram $\lambda^{(j)}$
is higher than that of $\mu^{(j)}$ for all $j$.

Now, for each path $p\in {\cal P}(\Lambda)$, 
there is a unique cylindrical preimage 
$\blambda \in \pi^{-1}(p)$ which is
higher than all other cylindrical preimages of $p$ \cite{DJKMO}. 
It is called the highest-lift of $p$ and is denoted by $\blambda(p)$.
The following result highlights the correspondence between paths
and their highest-lifts.

\begin{lemma}\label{WtLemma}{\rm\cite{DJKMO}}
If $p\in {\cal P}(\Lambda)$ then
$$
{\rm wt}(\blambda(p)) = {\rm wt}(p)\,.
$$
\end{lemma}
\Proof We outline the proof of this result in Appendix A. \cqfd

Let ${\cal Y}(\Lambda)$ denote the set of highest-lifts 
of the paths $p\in {\cal P}(\Lambda)$.
For example, a few paths $p\in\P(2\La_0)$ for $n=2$, and their
corresponding highest-lifts $\blambda \in {\cal Y}(2\La_0)$ are shown
in Figure~\ref{FIG0}.
Using the above result, we are now able to rewrite
Theorem~\ref{TH2.2} in the following form (also due to \cite{DJKMO}):
$$
{\rm ch} V(\Lambda) =
    \sum_{\blambda\in {\cal Y}(\Lambda)} e^{{\rm wt}(\blambda)} \,.
$$

The significance of paths and highest-lifts here is in fact explained by
the following theorem of Jimbo, Misra, Miwa and Okado.
\begin{figure}[t]
\begin{center}
\leavevmode
\epsfxsize =13.5cm
\epsffile{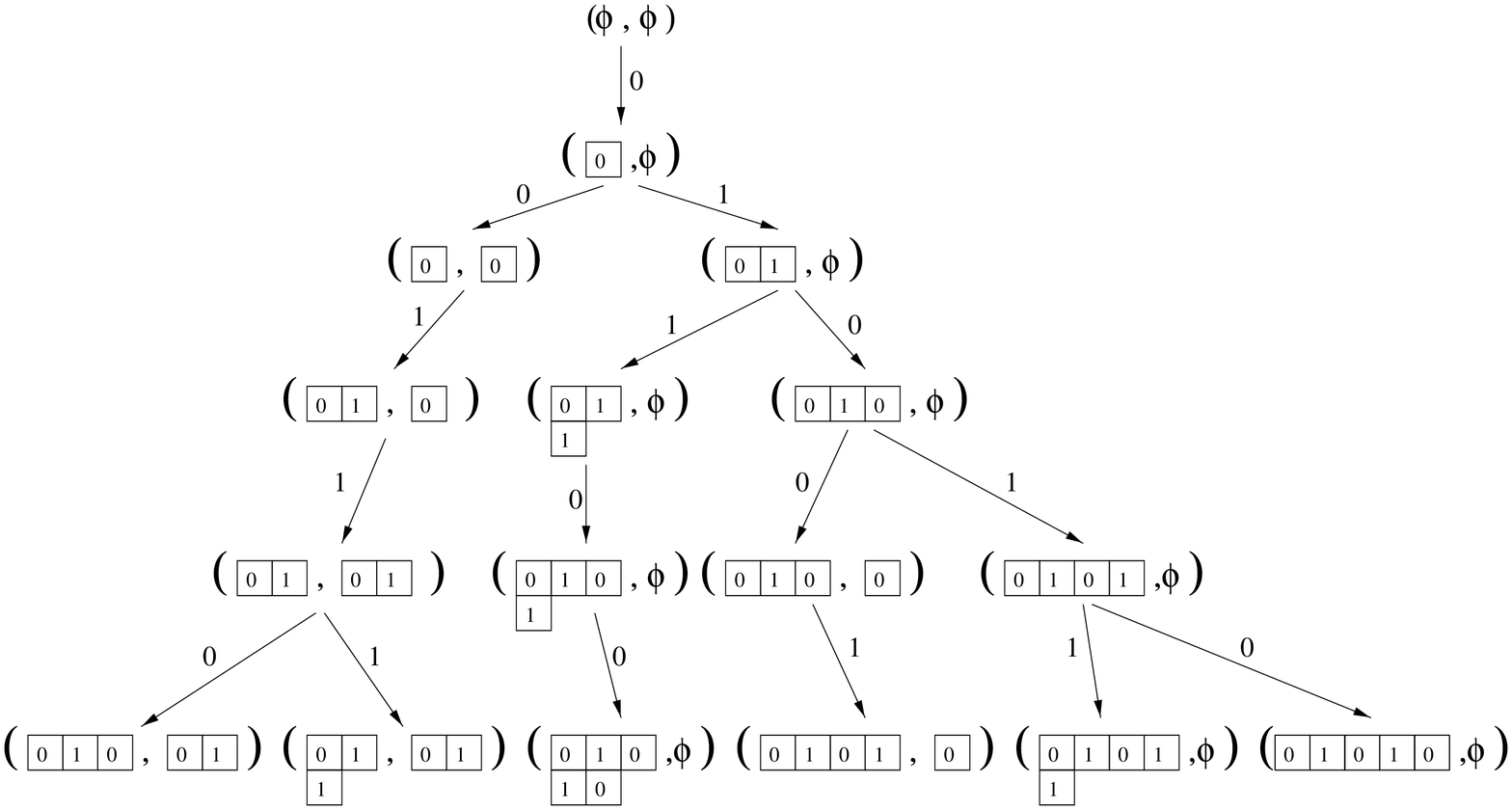}
\end{center}
\caption{\label{FIG1} The crystal graph the $U_q(\slchap_2)$-module
$V_q(2\La_0)$ labelled by ${\cal Y}(2\La_0)$}
\end{figure}
\begin{theorem}{\rm\cite{JMMO}}
The crystal graph of $V_q(\Lambda)$ is the full subgraph of
the crystal graph of $\F_q(\Lambda)$ with set of vertices
${\cal Y}(\Lambda)$.
\end{theorem}
Thus, denoting by $(L(\Lambda),B(\Lambda))$ the crystal basis
at $q=0$ of $V_q(\Lambda)$, we can identify $B(\Lambda)$ with
${\cal P}(\Lambda)$ or ${\cal Y}(\Lambda)$. 
The crystal graph of the $U_q(\slchap_2)$-module $V_q(2\La_0)$
with vertices labelled by ${\cal Y}(2\Lambda_0)$
is shown in Figure~\ref{FIG1}, up to the principal degree 5.

\medskip
We end this section by stating several useful results concerning
$\La$-paths and highest-lift multipartitions.
We first give a direct characterisation of highest-lift multipartitions.
It is shown in Appendix~A that this characterisation is
equivalent to that of \cite{JMMO}.

\begin{proposition}\label{TREVOR}
Let $\blambda$ be a cylindrical multipartition
of highest weight $\Lambda$. 
Then $\blambda \in {\cal Y}(\Lambda)$ if and only if
for all $k>0$, among the colours appearing at the right ends of
the length $k$ rows of $\blambda$,  at least one element
of $\{0,1,\ldots,n-1\}$ does not occur.
\end{proposition}
Note that when $\La=\La_i$ is a fundamental weight, we recover the  
usual characterisation of ${\cal Y}(\Lambda_i)$ as the
set of {\em $n$-regular partitions}, \ie partitions 
$\lambda = (\cdots 2^{m_2}1^{m_1})$ with all $m_i<n$.
\begin{example}\label{EXHL}{\rm
Take $n=4$, $\Lambda = 2\Lambda_0+\Lambda_1$, and let 
$\blambda$ be the multipartition with the following
coloured diagram:
$$
\smyoungd{\omit\phantom{\hrulefill}\cr
 \omit&\cr
 \multispan{22}\hrulefill\cr
 &0&&1&&2&&3&&0&&1&&2&&3&&0&&1&&\cr
 \multispan{21}\hrulefill\cr
 &3&&0&&1&&2&&3&&0&&1&&2&&3&&0&\cr
 \multispan{21}\hrulefill\cr
 &2&&3&&0&&1&&2&&3&&0&&1&\cr
 \multispan{17}\hrulefill\cr
 &1&&2&&3&&0&\cr
 \multispan{9}\hrulefill\cr
 &0&&1&&2&&3&\cr
 \multispan{9}\hrulefill\cr
 &\cr
 \multispan{1}\hrulefill\cr
 }
\quad
\smyoungd{\omit\phantom{\hrulefill}\cr
 \omit&\cr
 \multispan{22}\hrulefill\cr
 &0&&1&&2&&3&&0&&1&&2&&3&&0&&&\omit&\cr
 \multispan{19}\hrulefill\cr
 &3&&0&&1&&2&&3&&0&&1&&2&&3&\cr
 \multispan{19}\hrulefill\cr
 &2&\cr
 \multispan{3}\hrulefill\cr
 &1&\cr
 \multispan{3}\hrulefill\cr
 &\cr
 \multispan{1}\hrulefill\cr
 &\cr
 \multispan{1}\hrulefill\cr
 }
\quad
\smyoungd{\multispan{22}\hrulefill\cr
 &1&&2&&3&&0&&1&&2&&3&&0&&1&&2&&\cr
 \multispan{21}\hrulefill\cr
 &0&&1&&2&&3&&0&&1&&2&\cr
 \multispan{15}\hrulefill\cr
 &3&\cr
 \multispan{3}\hrulefill\cr
 &\cr
 \multispan{1}\hrulefill\cr
 &\cr
 \multispan{1}\hrulefill\cr
 &\cr
 \multispan{1}\hrulefill\cr
 &\cr
 \multispan{1}\hrulefill\cr
 }
$$
We immediately see that this is a cylindrical
multipartition. To check that it is a highest-lift
multipartition, we first note that 0 does not occur at the end of
a row of length 1. There are no rows of length 2, 3, 5, 6, and $k$ for
$k>10$ and so nothing to check in these cases.
For rows of length 4, 7, 8, 9, and 10, we see that, for
example, 1, 0, 0, 1 and 3 respectively do not occur.
Thus our example is also highest-lift.\endex}
\end{example}
The following relations between a $\Lambda$-path 
and its highest-lift are also proved in Appendix~A.
We introduce some notation.
Let $\sigma$ denote the linear map defined by
$\sigma(\La_i)=\La_{i-1}$.
For $i\in \Z$, put $\alpha'_i = -\Lambda_{i-1}+2\Lambda_i-\La_{i+1}$.
Note that $\sum_{i=0}^{n-1}\alpha'_i = 0$.
\begin{proposition}\label{EndColLemma}
Let $p\in{\cal P}(\Lambda)$ and $\blambda = \blambda(p)$. 
\begin{quote}
{\rm 1.} Let $\psi_r(k)$ be the number of $r$-nodes
at the end of length $k$ rows of $\blambda$.
Then,
\begin{equation}\label{det_colours}
\eta_{k-1}-\sigma(\eta_k)=
\sum_{j=0}^{l-1} \epsilon_{\gamma_j(k-1)}
-\sum_{j=0}^{l-1} \epsilon_{\gamma_j(k)-1}
=\sum_{r=0}^{n-1} \psi_{r}(k)\alpha^\prime_r.
\end{equation}
Since $\blambda$ is a highest-lift, one of the integers
$\psi_0(k),\ldots ,\psi_{n-1}(k)$ is $0$, and hence these
integers are unambiguously determined by Eq.~{\rm(\ref{det_colours})}.

{\rm 2.} Let $m_r(k)$ be the number of $r$-nodes of $\blambda$ 
belonging to columns to the right of the $k$th column
(including the $k$th column). Then
$$
p_{k-1}-\overline p_{k-1} = -\sum_{r=0}^{n-1} m_r(k) \alpha_r' \,.
$$
\end{quote}
\end{proposition}
\begin{example}{\rm
We continue Example~\ref{EXHL}. Let $p=\pi(\blambda)$ and
take $k=9$. We have
$$
\eta_8 - \sigma(\eta_9) = 
(\epsilon_0+2\epsilon_2)
- (2\epsilon_0 + \epsilon_2)
=-\epsilon_0+\epsilon_2=\alpha'_0+\alpha'_3.
$$
Correspondingly, there are two rows of length 9 having right end
nodes of colours 0 and 3 respectively.

We have $p_8 =  \La_1 + 2\La_2$ and $\overline p_8 = 2\La_0+\La_1$,
so that 
$$
p_8 - \overline p_8 = -2\La_0+2\La_2
=-\alpha'_0 + \alpha'_2
=-(3\alpha'_0 + 2\alpha'_1+\alpha'_2 +2\alpha'_3) \,.
$$
On the other hand, we can check that the 9th and 10th columns of $\blambda$
contain 8 nodes with colours $0,0,0,1,1,2,3,3$.
\endex}
\end{example}

Finally, we recall how to use paths or highest-lifts
to compute the length of a string of the
crystal graph of $V_q(\La)$. 
Let $b\in B(\Lambda)$ and let $\varepsilon_r(b)$
(\resp $\varphi_r(b)$) denote the greatest integer
$k$ such that $\tilde e_r^k(b) \not = 0$
(\resp $\tilde f_r^k(b) \not = 0$).
When $b$ is identified with $\blambda \in {\cal Y}(\Lambda)$
or $p\in {\cal P}(\Lambda)$, we write as well
$\varepsilon_r(b) = \varepsilon_r(\blambda)= \varepsilon_r(p)$
and $\varphi_r(b) = \varphi_r(\blambda)= \varphi_r(p)$.
Given $\nu = \sum_j \nu_j \epsilon_j\in {\cal A}_l^+$, we put
$\hat \nu = \sum_j \nu_j \Lambda_j\in P_l^+$.
For $A=\sum_j a_j \Lambda_j$ we define
$$
|A|^-_r =  \cases{-a_r & if $a_r\le 0$,\cr
0 & if $a_r >0 $.}
$$
\begin{proposition}{\rm\cite{JMMO}}\label{EPSILON}
\begin{quote}
{\rm 1.} For $\blambda \in {\cal Y}(\Lambda)$, $\varepsilon_r(\blambda)$
is equal to the number of normal removable $r$-nodes of~$\blambda$.

{\rm 2.} Let $p\in {\cal P}(\Lambda)$ and $\eta_k = p_{k+1}-p_k$
for all $k\ge0$.
Then 
\begin{equation}\label{EQ.EPS.PATH}
\varepsilon_r(p) = \max_{k\ge 0}\, |p_k-\hat\eta_k|_r^-\,.
\end{equation}
\end{quote}
\end{proposition} 
\begin{example}\label{EX2-15}{\rm
Let $n=3$, $\Lambda = 2\Lambda_1+\Lambda_2$.
Take $\blambda=((4,2),(3,1),(5))$, whose coloured
diagram is
$$
\smyoungd{\omit\phantom{\hrulefill}\cr
 \omit&\cr
 \multispan{10}\hrulefill\cr
 &1&&2&&0&&1&&\cr
 \multispan{9}\hrulefill\cr
 &0&&1&&\cr
 \multispan{5}\hrulefill\cr
 &\cr
 \multispan{1}\hrulefill\cr}
\quad
\smyoungd{\omit\phantom{\hrulefill}\cr
 \omit&\cr
 \multispan{8}\hrulefill\cr
 &1&&2&&0&&\cr
 \multispan{7}\hrulefill\cr
  &0&&\cr
 \multispan{3}\hrulefill\cr
&\cr
 \multispan{1}\hrulefill\cr}
\quad
\smyoungd{\multispan{12}\hrulefill\cr
 &2&&0&&1&&2&&0&&\cr
 \multispan{11}\hrulefill\cr
 &\cr
 \multispan{1}\hrulefill\cr
&\cr
 \multispan{1}\hrulefill\cr
&\cr
 \multispan{1}\hrulefill\cr}
$$
It is immediate from Proposition~{\ref{TREVOR}} that
$\blambda \in {\cal Y}(\Lambda)$.
The $r$-signatures of $\blambda$ are
\begin{eqnarray*}
&&r=0: \qquad R_{0,1}R_{3,1}R_{6,2}\\[1mm]
&&r=1: \qquad A_{1,2}A_{1,1}(R_{1,0}A_{4,1})(R_{4,0}A_{7,2})\\[1mm]
&&r=2: \qquad A_{-1,1}A_{-1,0}A_{2,0}A_{5,0}
\end{eqnarray*}
Hence $\varepsilon_0(\blambda) = 3$ and
$\varepsilon_1(\blambda) = \varepsilon_2(\blambda) =0$.
On the other hand, the path $p=\pi(\blambda)$ is given by
\begin{eqnarray*}
&&\eta=(122,012,022,011,222,001,112,022,001,\ldots ), \\[1mm]
&&p = (-3\Lambda_0+2\Lambda_1+4\Lambda_2,
-\Lambda_0+\Lambda_1+3\Lambda_2,
-\Lambda_0+\Lambda_1+3\Lambda_2,
2\Lambda_1+\Lambda_2, \\
&& \quad -\Lambda_0+\Lambda_1+3\Lambda_2, 
2\Lambda_0+\Lambda_1,
2\Lambda_1+\Lambda_2,
\Lambda_0+2\Lambda_2,
2\Lambda_0+\Lambda_1, \ldots ), \\[1mm]
&&p-\hat\eta = (-3\Lambda_0+\Lambda_1+2\Lambda_2,
-2\Lambda_0+2\Lambda_2, 
-2\Lambda_0+\Lambda_1+\Lambda_2, 
-\Lambda_0+\Lambda_2, \\
&& \quad -\Lambda_0+\Lambda_1, 
0,0, 0, 0, \ldots ).
\end{eqnarray*}
Hence, Proposition \ref{EPSILON}(2) gives 
$$\varepsilon_0(p) = \max(3,2,2,1,1,0,0,\ldots ) = 3\,,\qquad
\varepsilon_1(p) = \varepsilon_2(p) =0\,.$$
\endex}
\end{example}

\subsection{Tensor products, restricted paths and restricted 
highest-lifts}
\label{SECTTP}

The tensor product of two $U_q(\slchap_n)$-modules carries the structure
of a $U_q(\slchap_n)$-module defined via the coproduct
\begin{equation}\label{COM}
\Delta(q^h)=q^h\otimes q^h,\quad 
\Delta(e_i)=e_i\otimes q^{-h_i} + 1\otimes e_i, \quad
\Delta(f_i)=f_i\otimes 1 + q^{h_i}\otimes f_i.
\end{equation}
Let $(L_1,B_1)$ and $(L_2,B_2)$ be crystal bases of the integrable
$U_q(\slchap_n)$-modules $M_1$ and $M_2$.
Let $B_1\otimes B_2$ denote the basis $\{ u\otimes v,\ u\in B_1,\,
v\in B_2\}$
of $(L_1/qL_1 )\otimes (L_2/qL_2)$. 
Then, Kashiwara \cite{Ka1,Ka2} has proved that
$(L_1\otimes L_2,\,B_1\otimes B_2)$ is a
crystal basis of $M_1\otimes M_2$, with the action of
$\tilde e_i$, $\tilde f_i$ on $B_1\otimes B_2$ given by
\begin{eqnarray}
\tilde{f_i}(u\otimes v)&=&\cases{
\tilde{f_i}u\otimes v&if $\varphi_i(u)>\varepsilon_i(v)$,\cr
u\otimes\tilde{f_i}v&otherwise,\cr}\\
\tilde{e_i}(u\otimes v)&=&\cases{
\tilde{e_i}u\otimes v&if $\varphi_i(u)\ge\varepsilon_i(v)$,\cr
u\otimes\tilde{e_i}v&otherwise.\cr} \label{ETILD}
\end{eqnarray}
This gives a powerful way of computing tensor product
multiplicities for $U_q(\slchap_n)$ and $\slchap_n$.
Let $\Lambda^\prime, \Lambda^{\prime\prime}, \Lambda \in P^+$.
The multiplicity
$c_{\Lambda^\prime\,\Lambda^{\prime\prime}}^\Lambda$
of $V(\Lambda)$ in the
tensor product $V(\Lambda^\prime)\otimes V(\Lambda^{\prime\prime})$
is equal to the number of vertices $b_1\otimes b_2$ of
$B(\Lambda^\prime)\otimes B(\Lambda^{\prime\prime})$ that satisfy
$$
\wt(b_1\otimes b_2) = \Lambda
\quad \mbox{and} \quad \tilde{e}_i (b_1\otimes b_2) = 0, 
\qquad (i=0,\ldots ,n-1).
$$
By (\ref{ETILD}), this last condition is equivalent to the
fact that $b_1=b_{\Lambda^\prime}$,  the origin of the crystal graph
of $V(\Lambda^\prime)$, and
\begin{equation}
\varepsilon_i(b_2)\le \varphi_i(b_{\Lambda^\prime}) = 
\<\Lambda^\prime,h_i\> , \qquad (i=0,1,\ldots,n-1).
\end{equation}
Hence we get,
\begin{proposition}\label{TPcorollary}
The multiplicity $c_{\Lambda^\prime\,\Lambda^{\prime\prime}}^\Lambda$
of $V(\Lambda)$ in
$V(\Lambda^\prime)\otimes V(\Lambda^{\prime\prime})$
is equal to the number of vertices $b_2$ of $B(\Lambda^{\prime\prime})$
such that
\begin{equation}\label{TPMULT}
\wt(b_2)=\Lambda - \Lambda^\prime \quad \mbox{  and  } \quad
\varepsilon_i(b_2)\le \<\Lambda^\prime,h_i\> , \qquad (i=0,1,\ldots,n-1).
\end{equation}
\end{proposition}

{}From now on, we fix positive integers $l^{\prime}, l^{\prime\prime}$,
$l=l^{\prime}+ l^{\prime\prime}$, and weights $\LP\in P^+_{l^\prime}$,
$\LPP\in P^+_{l^{\prime\prime}}$ and $\Lambda \in P^+_l$.
We define the {\em branching function}
\begin{equation}\label{BRFU}
b^\Lambda_{\Lambda'\Lambda''}(z) = 
\sum_k z^k c^{\Lambda-k\delta}_{\Lambda'\Lambda''}\,.
\end{equation}
Let us identify the crystal basis $B(\Lambda^{\prime\prime})$ with the set 
${\cal P}(\Lambda^{\prime\prime})$ of $\Lambda^{\prime\prime}$-paths,
as explained in Section~\ref{SECT2.1}.
Using Proposition~\ref{TPcorollary} and Eq.~(\ref{EQ.EPS.PATH}) of
Proposition~\ref{EPSILON}, it was proved in \cite{JMMO} that 
$b^\Lambda_{\Lambda'\Lambda''}(z)$ is 
equal to the generating function (with respect to the energy $E(p)$)
of certain subsets of ${\cal P}(\Lambda^{\prime\prime})$
coming from a class of restricted-solid-on-solid solvable models
(RSOS models). We proceed to describe those paths in detail.

\medskip
A level $l'+l''$ path $p=(p_k)$ is said to be 
{\em $(\LP,\LPP)$-restricted} if
\begin{enumerate}
\item $p-\LP:=(p_k-\LP) \in {\cal P}(\LPP)$, which
implies that 
$p_{k+1}-p_k = \eta_k \in {\cal A}_{l''}^+$; 
\item $p_k-\hat\eta_k\in P^+_{l^\prime}$ for all $k\ge0$.
\end{enumerate}
The set of $(\LP,\LPP)$-restricted paths is denoted by 
${\cal P}(\LP,\LPP)$.
The ground state path $\overline p\in{\cal P}(\LP,\LPP)$
is defined to be $\overline p=\overline p^{\prime\prime}+\LP$
where $\overline p^{\prime\prime}$ is the ground state path of
${\cal P}(\LPP)$.
\begin{example}\rm
We consider again the path $p\in\P(2\La_1+\La_2)$ of Example~\ref{EX2-15}.
The path $p+\LP$ is an element of ${\cal P}(\LP,2\La_1+\La_2)$
for all $\LP = a_0\La_0+a_1\La_1+a_2\La_2$ with
$a_0 \ge 3 , \ a_1\ge 0, \ a_2\ge 0.$
\endex
\end{example}
Note that if $p\in{\cal P}(\LP,\LPP)$,
then $p_k\in P^+_{l^\prime+l^{\prime\prime}}$ for all $k\ge0$.
Note also that in the case where $\LPP=\Lambda_0$,
the previous definition is equivalent to Definition 2.12.1 
of \cite{FLOTW2}.

Since $\{p-\LP\mid p\in{\cal P}(\LP,\LPP)\}\subset{\cal P}(\LPP)$,
we can associate to the restricted path $p\in {\cal P}(\LP,\LPP)$
the highest-lift multipartition $\blambda(p-\LP)$, that we
also denote by $\blambda(p)$, with a slight abuse of notation.
Similarly, write $\ell(p) = \ell(p-\LP)$ and call it the 
{\em length} of $p\in \P(\LP,\LPP)$.
Finally let
$$
{\cal Y}(\LP,\LPP)
=\left\{\blambda(p)\mid p\in{\cal P}(\LP,\LPP)\right\}
$$
denote the set of 
{\em $(\LP,\LPP)$-restricted highest-lift multipartitions}.
Equivalently, 
$$
{\cal Y}(\LP,\LPP)
=\left\{\blambda \in {\cal Y}(\LPP) \mid 
\varepsilon_i(\blambda) \le \langle \LP,h_i\rangle ,
\ 0\le i \le n-1 \right\}\,.
$$
\begin{theorem}{\rm\cite{JMMO}}
The branching function is given by
$$
b^{\Lambda}_{\Lambda^\prime,\Lambda^{\prime\prime}}(z)
=\sum_{
\scriptstyle p\in{\cal P}(\Lambda^\prime,\Lambda^{\prime\prime})\atop
\scriptstyle p_0=\Lambda}
z^{E(p)}
=\sum_{
\scriptstyle\sblambda\in{\cal Y}(\Lambda^\prime,\Lambda^{\prime\prime})\atop
\scriptstyle\wt(\blambda)=\Lambda -\Lambda^\prime\,(\mod\delta)}
z^{N^0(\sblambda)}.
$$
\end{theorem}
\begin{example}
\rm
Take $n=2$ and consider the branching function
$b^{3\La_0}_{\La_0,2\La_0}(z)$.
As can be seen from Figure~\ref{FIG1}, the first elements
of ${\cal Y}(\La_0,2\La_0)$ that contribute to this function
are 
$$
(\emptyset \,,\, \emptyset)\,, \qquad
((3,1)\,,\,\emptyset)\,,
$$
whence 
$b^{3\La_0}_{\La_0,2\La_0}(z) = 1 + z^2 + \mbox{higher terms}$.
Similarly, the first elements 
of ${\cal Y}(\La_0,2\La_0)$ that contribute to 
$b^{\La_0+2\La_1}_{\La_0,2\La_0}(z)$ are
$$
((1)\,,\,\emptyset)\,,\qquad
((3)\,,\,\emptyset)\,,\qquad
((3)\,,\,(2))\,,\qquad
((5)\,,\,\emptyset)\,,\qquad
$$
which gives
$b^{\La_0+2\La_1}_{\La_0,2\La_0}(z) = z + z^2 + 2z^3 +
\mbox{higher terms}$.

Of course, one may use paths instead of multipartitions to
obtain the same results.
\endex
\end{example}
When $\LP$ or $\LPP$ is a fundamental weight, there is also
a closed formula for $b^{\Lambda}_{\Lambda^\prime,\Lambda^{\prime\prime}}(z)$
in terms of theta-functions (see below, Section~\ref{SECT3.3}).

\subsection{ A bijection between ${\cal Y}(\LP,\LPP)$ 
and ${\cal Y}(\sharp\LPP,\sharp\LP)$}\label{SECT2.3}

As indicated in Section~\ref{SECTTP}, the tensor product
$V(\LP)\otimes V(\LPP)$ may be calculated by enumerating
certain elements of $\P(\LPP)$ or ${\cal Y}(\LPP)$.
On the other hand, since
$V(\LP)\otimes V(\LPP)$ is isomorphic to $V(\LPP)\otimes V(\LP)$,
this tensor product may equally be calculated using 
$\P(\LP)$ or ${\cal Y}(\LP)$.
Thus there exist weight-preserving bijections between
$\P(\LP,\LPP)$ and $\P(\LPP,\LP)$ or between 
${\cal Y}(\LP,\LPP)$ and ${\cal Y}(\LPP,\LP)$.
It turns out that such bijections are not easy to describe in terms
of the paths or the partitions involved.
Moreover, one can check from small examples
that such bijections cannot preserve the length of paths.

However, there exists an equally useful but much more easily
described length-preser\-ving bijection
between $\P(\LP,\LPP)$ and $\P(\sharp\LPP,\sharp\LP)$,
where $\sharp$ is the root diagram automorphism exchanging
$\Lambda_i$ and $\Lambda_{-i}$.
This bijection is the manifestation of a simple correspondence between
the tensor products of Fock spaces $\F_q(\LP)\otimes \F_q(\LPP)$
and $\F_q(\sharp\LPP)\otimes \F_q(\sharp\LP)$,
which we discuss in Appendix B.

\begin{definition}\nobreak
Let 
$\Lambda=\sum_{i=0}^{n-1} a_i\Lambda_i \in P_l$,
and $p$ be a level $l$ path,
with
$p_k=\sum_{i=0}^{n-1} a_i(k)\Lambda_i$.
We put
$\sharp\Lambda=\sum_{i=0}^{n-1} a_i\Lambda_{-i}$, and
we define the path $\sharp p$ by
$(\sharp p)_k=\sum_{i=0}^{n-1} a_i(k)\Lambda_{k-i}.$
\end{definition}
Note that $(\sharp p)_k$  differs from $\sharp (p_k)$
by a shift of indices. 
Hereafter, $\sharp p_k$ will always mean $(\sharp p)_k$.

\begin{proposition}\label{image_lemma}
Let $\LP\in P_{l'}^+$ and $\LPP\in P_{l''}$.
The map $p\mapsto \sharp p$ defines a bijection between
${\cal P}(\LP,\LPP)$ and ${\cal P}(\sharp\LPP,\sharp\LP)$.
Furthermore, the paths $p$ and $\sharp p$ have the same length.
\end{proposition}

\Proof For fixed $k$, we have
$$
p_{k}=
\Lambda_{w_0(k)}+\Lambda_{w_1(k)}+\cdots+\Lambda_{w_{l'+l''-1}(k)},
$$
where $0\le w_i(k)<n$.
Since $p$ is a $(\LP,\LPP)$-restricted path,
we may assume that 
$w_j(k)=\gamma_{j-l^\prime}(k)$ for $l^\prime\le j<l'+l''$. Then
$$
p_{k+1}=
\Lambda_{w_0(k)}+\cdots+\Lambda_{w_{l^\prime-1}(k)}
+\Lambda_{w_{l^\prime}(k)+1}+\cdots+\Lambda_{w_{l'+l''-1}(k)+1},
$$
so that
\[
\mathalign{
&\sharp\eta_k=(\sharp p_{k+1}-\sharp\LPP)-(\sharp p_k-\sharp\LPP)
=\sharp p_{k+1}-\sharp p_k\cr
&\qquad=(\Lambda_{k+1-w_0(k)}+\cdots+\Lambda_{k+1-w_{l^\prime-1}(k)}
+\Lambda_{k-w_{l^\prime}(k)}+\cdots+\Lambda_{k-w_{l'+l''\!-1}(k)})\cr
&\qquad\quad
- (\Lambda_{k-w_0(k)}+\cdots+\Lambda_{k-w_{l^\prime-1}(k)}
+\Lambda_{k-w_{l^\prime}(k)}+\cdots+\Lambda_{k-w_{l'+l''-1}(k)})\cr
&\qquad=
(\Lambda_{k+1-w_0(k)}+\cdots+\Lambda_{k+1-w_{l^\prime-1}(k)})
- (\Lambda_{k-w_0(k)}+\cdots+\Lambda_{k-w_{l^\prime-1}(k)})\cr
&\qquad=
\epsilon_{k-w_0(k)}+\cdots+\epsilon_{k-w_{l^\prime-1}(k)}\,. 
\cr}
\]
Thus $\sharp\eta_k\in{\cal A}_{l'}^+$.
Furthermore, on setting $\LP=\La_{u_0}+\cdots+\La_{u_{l^\prime-1}}$
and $\LPP=\La_{v_0}+\cdots+\La_{v_{l''-1}}$, we have
for each $k\ge \ell(p)$,
$$
p_k=\Lambda_{u_0}+\cdots+\Lambda_{u_{l^\prime-1}}
+\Lambda_{v_0+k}+\cdots+\Lambda_{v_{l''-1}+k},
$$
whence
$$
\sharp p_k=\Lambda_{k-u_0}+\cdots+\Lambda_{k-u_{l^\prime-1}}
+\Lambda_{-v_0}+\cdots+\Lambda_{-v_{l''-1}},
$$
which shows that $\sharp p-\sharp\LPP$ is an element of 
${\cal P}(\sharp\LP)$.
It also shows that $\ell(\sharp p-\sharp\LPP)\le \ell(p-\LP)$.
Exchanging the roles of
$p$ and $\sharp p$ yields $\ell(p-\LP)\le \ell(\sharp p-\sharp\LPP)$,
and therefore $\ell(p-\LP)= \ell(\sharp p-\sharp\LPP)$.

Now, comparing the expression 
$$
\sharp p_k=\La_{k-w_0(k)}+\cdots+\La_{k-w_{l'+l''-1}(k)}
$$
with the above expression for $\sharp \eta_k$, we see that $\sharp p$ is 
$(\sharp\LPP,\sharp\LP)$-restricted.
Finally the bijectivity of $p\mapsto \sharp p$ follows
from the fact that it is an involution between two 
finite sets with the same cardinality. \cqfd

We now map the paths $p$ and $\sharp p$ to their corresponding
multipartitions, thereby establishing a bijection between
${\cal Y}(\LP,\LPP)$ and ${\cal Y}(\sharp\LPP,\sharp\LP)$.
Let $\blambda=\blambda(p)$ be the highest-lift 
of $p\in{\cal P}(\LP,\LPP)$
and define $\sharp\blambda=\blambda(\sharp p)$ to be the
highest-lift of $\sharp p$.
\begin{theorem}\label{diag_bij}
$\sharp\blambda$ is the unique element of ${\cal Y}(\sharp\LPP,\sharp\LP)$ 
for which, for all $k>0$ and $0\le r<n$,
the number of length $k$ rows with right end of colour $r$ 
is equal to the number of length $k$ rows of $\blambda$
with left end of colour $-r$.
In particular, $\blambda$ and $\sharp\blambda$ have the same
number of rows of length $k$, for all $k>0$.
\end{theorem}

\Proof 
Let $k>0$ and write $p_{k-1}$ in the form 
$$
p_{k-1}=\sum_{j=0}^{l'+l''-1} \Lambda_{w_j(k-1)}.
$$
Then
$$
\mathalign{
p_{k}&=\sum_{j=0}^{l'+l''-1} \Lambda_{w_j(k-1)}
+\sum_{j=0}^{l''-1}\epsilon_{\gamma_j(k-1)},\cr
p_{k+1}&=\sum_{j=0}^{l'+l''-1} \Lambda_{w_j(k-1)}
+\sum_{j=0}^{l''-1}\epsilon_{\gamma_j(k-1)}
+\sum_{j=0}^{l''-1}\epsilon_{\gamma_j(k)}.\cr}
$$
 Applying the map $p\mapsto \sharp p$ yields
(taking care to transform the $\epsilon$ terms correctly)
$$
\mathalign{
\sharp p_{k-1}&=\sum_{j=0}^{l'+l''-1} \Lambda_{k-1-w_j(k-1)},\cr
\sharp p_{k}&=\sum_{j=0}^{l'+l''-1} \Lambda_{k-w_j(k-1)}
-\sum_{j=0}^{l''-1}\epsilon_{k-1-\gamma_j(k-1)},\cr
\sharp p_{k+1}&=\sum_{j=0}^{l'+l''-1} \Lambda_{k+1-w_j(k-1)}
-\sum_{j=0}^{l''-1}\epsilon_{k-\gamma_j(k-1)}
-\sum_{j=0}^{l''-1}\epsilon_{k-\gamma_j(k)}.\cr}
$$
Therefore
$$
\mathalign{
\sharp\eta_k&=
\sharp p_{k+1}-\sharp p_k=\sum_{j=0}^{l'+l''-1} \epsilon_{k-w_j(k-1)}
+\sum_{j=0}^{l''-1}(\epsilon_{k-1-\gamma_j(k-1)}
-\epsilon_{k-\gamma_j(k-1)}
-\epsilon_{k-\gamma_j(k)})\cr
\noalign{\noindent and}
\sharp\eta_{k-1}&=
\sharp p_k-\sharp p_{k-1}=\sum_{j=0}^{l'+l''-1} \epsilon_{k-1-w_j(k-1)}
-\sum_{j=0}^{l''-1}\epsilon_{k-1-\gamma_j(k-1)},\cr}
$$
whereupon, if $\psi_r^\sharp(k)$ is the number of $r$-nodes 
at the right ends of the length $k$ rows of $\sharp\blambda$,
on using Proposition~\ref{EndColLemma} Eq.~(\ref{det_colours}), we obtain:
$$
\sum_{r=0}^{n-1} \psi^\sharp_r(k) \alpha'_r=
\sharp \eta_{k-1} - \sigma(\sharp\eta_k) =
-\sum_{j=0}^{l''-1}\epsilon_{k-2-\gamma_j(k-1)}
+\sum_{j=0}^{l''-1}\epsilon_{k-1-\gamma_j(k)}.
$$
In the case of $\blambda$, expanding (\ref{det_colours}) gives:
$$
 \sum_{j=0}^{l''-1}\epsilon_{\gamma_j(k-1)}
-\sum_{j=0}^{l''-1}\epsilon_{\gamma_j(k)-1}
=\sum_{r=0}^{n-1} \psi_r(k) \alpha'_r.
$$
On applying the linear map $\La_i \mapsto \La_{k-1-i}$,
so that $\epsilon_i \mapsto -\epsilon_{k-2-i}$
and $\alpha_i \mapsto \alpha_{k-1-i}$, 
to both sides of this formula, we obtain:
$$
-\sum_{j=0}^{l''-1} \epsilon_{k-2-\gamma_j(k-1)}
+\sum_{j=0}^{l''-1} \epsilon_{k-1-\gamma_j(k)}
=\sum_{r=0}^{n-1} \psi_r(k)\alpha^\prime_{k-1-r}.
$$
Thus $\psi_r^\sharp(k)=\psi_{k-1-r}(k)$ for $0\le r<n$.
Since $\psi_{k-1-r}(k)$ is the multiplicity of the
colour $-r$ at the left ends of the length $k$ rows of $\blambda$,
this multiplicity is equal to the multiplicity of the colour $r$
at the right ends of the length $k$ rows of $\sharp \blambda$.

That $\sharp \blambda$ is determined uniquely follows since
the colour of the nodes at the right ends of the length
$k$ rows determine, via Eq.~(\ref{det_colours}),
the path $\pi(\sharp \blambda)$, which in turn 
determines $\sharp \blambda$ since it is a highest-lift multipartition.
\cqfd

\begin{example}{\rm
To illustrate Theorem \ref{diag_bij}, consider for $n=3$ the 
highest-lift multipartition 
$\blambda \in {\cal Y}(\Lambda_1+2\Lambda_2,\Lambda_0+\Lambda_1)$
whose coloured diagram is
$$
\smyoungd{\omit\phantom{\hrulefill}\cr
 \omit&\cr
 \multispan{20}\hrulefill\cr
 &0&&1&&2&&0&&1&&2&&0&&1&&2&&\cr
 \multispan{19}\hrulefill\cr
 &2&&0&&1&&2&&0&&1&&2&&0&\cr
 \multispan{17}\hrulefill\cr
 &1&&2&&0&&1&&2&&0&&1&\cr
 \multispan{15}\hrulefill\cr
 &0&&1&&2&&0&&1&\cr
 \multispan{11}\hrulefill\cr
 &2&&0&&1&&2&\cr
 \multispan{9}\hrulefill\cr
 &1&&2&&0&&1&\cr
 \multispan{9}\hrulefill\cr
 &0&\cr
 \multispan{3}\hrulefill\cr
 &2&\cr
 \multispan{3}\hrulefill\cr
 &\cr
 \multispan{1}\hrulefill\cr
 }
\quad
\smyoungd{\multispan{20}\hrulefill\cr
 &1&&2&&0&&1&&2&&0&&1&&2&&0&&\cr
 \multispan{19}\hrulefill\cr
 &0&&1&&2&&0&&1&&2&&0&&1&&2&&\cr
 \multispan{19}\hrulefill\cr
 &2&&0&&1&&2&&0&&1&&2&\cr
 \multispan{15}\hrulefill\cr
 &1&&2&&0&&1&&2&&0&\cr
 \multispan{13}\hrulefill\cr
 &0&&1&&2&&0&&1&\cr
 \multispan{11}\hrulefill\cr
 &2&&0&&1&\cr
 \multispan{7}\hrulefill\cr
 &1&&2&&0&\cr
 \multispan{7}\hrulefill\cr
 &\cr
 \multispan{1}\hrulefill\cr
 &\cr
 \multispan{1}\hrulefill\cr
 &\cr
 \multispan{1}\hrulefill\cr
 }
$$
By definition, $\sharp\blambda$ is a multipartition
of highest weight $2\Lambda_1+\Lambda_2$ comprising 3 partitions.
{}From Theorem \ref{diag_bij}, this multipartition must have
3 rows of length 9, with left end of colour
0, 1 and 1. These can only be placed as the first
rows of the 3 partitions in one way so as to obtain a
cylindrical multipartition: two in the first partition and one
in the second. Next, a row of length 8, beginning with colour
0 can only be placed in the second partition to obtain
a cylindrical multipartition. Proceeding in this way, distributing
the rows in order of decreasing length, we obtain the following
multipartition
$$
\smyoungd{\omit\phantom{\hrulefill}\cr
 \omit&\cr
 \multispan{20}\hrulefill\cr
 &1&&2&&0&&1&&2&&0&&1&&2&&0&&\cr
 \multispan{19}\hrulefill\cr
 &0&&1&&2&&0&&1&&2&&0&&1&&2&&\cr
 \multispan{19}\hrulefill\cr
 &2&&0&&1&&2&&0&\cr
 \multispan{11}\hrulefill\cr
 &1&&2&&0&&1&\cr
 \multispan{9}\hrulefill\cr
 &0&&1&&2&\cr
 \multispan{7}\hrulefill\cr
 &2&&0&&1&\cr
 \multispan{7}\hrulefill\cr
 &\cr
 \multispan{1}\hrulefill\cr}
\quad
\smyoungd{\omit\phantom{\hrulefill}\cr
 \omit&\cr
 \multispan{20}\hrulefill\cr
 &1&&2&&0&&1&&2&&0&&1&&2&&0&&\cr
 \multispan{19}\hrulefill\cr
 &0&&1&&2&&0&&1&&2&&0&&1&&\cr
 \multispan{17}\hrulefill\cr
 &2&&0&&1&&2&&0&\cr
 \multispan{11}\hrulefill\cr
 &1&\cr
 \multispan{3}\hrulefill\cr
 &0&\cr
 \multispan{3}\hrulefill\cr
 &\cr
 \multispan{1}\hrulefill\cr
 &\cr
 \multispan{1}\hrulefill\cr}
\quad
\smyoungd{\multispan{20}\hrulefill\cr
 &2&&0&&1&&2&&0&&1&&2&&&\omit&&\omit&\cr
 \multispan{15}\hrulefill\cr
 &1&&2&&0&&1&&2&&0&&1&\cr
 \multispan{15}\hrulefill\cr
 &0&&1&&2&&0&&1&&2&\cr
 \multispan{13}\hrulefill\cr
 &2&&0&&1&&2&\cr
 \multispan{9}\hrulefill\cr
 &\cr
 \multispan{1}\hrulefill\cr
 &\cr
 \multispan{1}\hrulefill\cr
 &\cr
 \multispan{1}\hrulefill\cr
 &\cr
 \multispan{1}\hrulefill\cr}
$$
\endex}
\end{example}
In fact, when $l'=1$ a simple
algorithm determines $\sharp\blambda$.
Indeed, $\blambda$ is an $l''$-tuple
multipartition, and from the second part of
Theorem \ref{diag_bij}, $\sharp\blambda$ is a single partition 
obtained by collecting together all the rows of
the constituent partitions of $\blambda$.
\begin{example}{\rm
Let $\blambda$ be the multipartition
of Example~\ref{EXHL}. 
We can check that $\blambda \in {\cal Y}(\La_0, 2\La_0+\La_1)$.
Then $\sharp\blambda \in {\cal Y}(2\La_0+\La_3, \La_0)$
is the following single partition
$$
\smyoungd{\multispan{22}\hrulefill\cr
 &0&&1&&2&&3&&0&&1&&2&&3&&0&&1&&\cr
 \multispan{21}\hrulefill\cr
 &3&&0&&1&&2&&3&&0&&1&&2&&3&&0&\cr
 \multispan{21}\hrulefill\cr
 &2&&3&&0&&1&&2&&3&&0&&1&&2&&3&\cr
 \multispan{21}\hrulefill\cr
 &1&&2&&3&&0&&1&&2&&3&&0&&1&\cr
 \multispan{19}\hrulefill\cr
 &0&&1&&2&&3&&0&&1&&2&&3&&0&\cr
 \multispan{19}\hrulefill\cr
 &3&&0&&1&&2&&3&&0&&1&&2&\cr
 \multispan{17}\hrulefill\cr
 &2&&3&&0&&1&&2&&3&&0&\cr
 \multispan{15}\hrulefill\cr
 &1&&2&&3&&0&\cr
 \multispan{9}\hrulefill\cr
 &0&&1&&2&&3&\cr
 \multispan{9}\hrulefill\cr
 &3&\cr
 \multispan{3}\hrulefill\cr
 &2&\cr
 \multispan{3}\hrulefill\cr
 &1&\cr
 \multispan{3}\hrulefill\cr
 &\cr
 \multispan{1}\hrulefill\cr
 }
$$
One can check that, for each $k$, 
the colours of the leftmost nodes of the rows of
length $k$ are the $n$-complements of
those of the rightmost nodes of the rows of length $k$ in $\blambda$.
\endex}\end{example}
%

\subsection{Another description of ${\cal Y}(\La_u,\La)$}
\label{SECT2.4}

Let $\La_u$ be a fixed fundamental weight and let $\La\in P_l^+$.
By means of Proposition~\ref{EPSILON}, the
multipartitions $\blambda \in {\cal Y}(\La)$ which belong
to ${\cal Y}(\La_u,\La)$, may be obtained by calculating the integers
$\varepsilon_r(\blambda)$ using the normal nodes of $\blambda$,
and checking that $\varepsilon_u(\blambda)\le 1$ and
$\varepsilon_j(\blambda)= 0$ for $j\not = u$.
In this section, we apply Theorem~\ref{diag_bij} to
obtain an alternative description of ${\cal Y}(\La_u,\La)$,
similar to that obtained in \cite{FOW} in the case where
$\La$ has level $l=1$.
\begin{theorem}\label{mult_char}
Let $\blambda\in{\cal Y}(\La)$
and let $\underline\lambda$ be the partition which, 
for each $k>0$, has as many rows of length $k$ as $\blambda$.
Then $\blambda\in{\cal Y}(\Lambda_{u},\La)$ if and only if
$\underline\lambda \in {\cal Y}(\La_{-u})$ and, for each $k>0$,
the number of length $k$ rows of $\underline\lambda$ with left end
of colour $r$
is equal to the number of length $k$ rows of $\blambda$
with right end of colour $-r$.
\end{theorem}
\Proof
The `only if' part follows immediately from Theorem \ref{diag_bij}.
Let us prove the converse.
Suppose that $\blambda\in{\cal Y}(\La)$ is such that
$\underline\lambda \in {\cal Y}(\La_{-u})$ and, for each $k>0$,
the number of length $k$ rows of $\underline\lambda$ with left end
of colour $r$
is equal to the number of length $k$ rows of $\blambda$
with right end of colour $-r$.
Put $p=\pi(\blambda) \in \P(\La)$ and write 
$$
p_k=\sum_{j=0}^{n-1}\beta_j(k)\La_j.
$$
We first prove the following
\begin{lemma}\label{LEMME}
The path $\underline p:= \pi(\underline\lambda) \in \P(\La_{-u})$ is given by
$$
\underline p= \overline p_{\La_{-u}} -\sharp\La + \sharp p \,,
$$
where $\overline p_{\La_{-u}}$ denotes the ground state
path of $\P(\La_{-u})$.
\end{lemma}
{\it Proof of Lemma}~\ref{LEMME}: 
Since $p_k-\hat\eta_k$ is not dominant for only a finite number of $k$
(it is always dominant if $k\ge k_*$, the length of $p$), we can
find a dominant integral weight $\chi$ such that
$p_k+\chi-\hat\eta_k$ is always dominant. Then
$p+\chi\in{\cal P}(\chi,\La)$ with the actual points given by
$$
(p+\chi)_k=\sum_{j=0}^{n-1}(\beta_j(k)+x_j)\La_j,
$$
where $\chi=\sum_{j=0}^{n-1} x_j\Lambda_j$.
Thence, by Proposition \ref{image_lemma},
$\sharp(p+\chi)\in{\cal P}(\sharp\La,\sharp\chi)$.
Theorem \ref{diag_bij} shows that the multipartition
$\widehat\blambda\in{\cal Y}(\sharp\La,\sharp\chi)$ corresponding
to this path has, for each $r$ and $k$, as many length $k$ rows 
with right end of
colour $r$ as does $\blambda$ with left end of colour $-r$.
The same multipartition corresponds to the $\sharp\chi$-path
$\sharp(p+\chi)-\sharp\La$.

Let $m_j(k)$ be the number of
$j$-nodes in and to the right of the $k$th column
of the partitions composing $\widehat\blambda$.
By Proposition~\ref{EndColLemma}, we have
$$
(\sharp(p+\chi)-\sharp\La)_k=
\sum_{j=0}^{n-1}x_j\Lambda_{k-j}
-\sum_{j=0}^{n-1}m_j(k+1)\alpha_j^\prime.
$$
{}From the assumption on $\underline\lambda$,
the rows of $\widehat\blambda$ and $\underline\lambda$ are
in bijection, and the colours of the nodes are preserved. 
Therefore, 
$$
\mathalign{
\underline p_k&=\Lambda_{k-u}-\sum_{j=0}^{n-1}m_j(k+1)\alpha_j^\prime\cr
&=\Lambda_{k-u}
+(\sharp(p+\chi)-\sharp\La)_k-\sum_{j=0}^{n-1}x_j\Lambda_{k-j}\cr
&=\Lambda_{k-u}-\sharp\La+\sum_{j=0}^{n-1}\beta_j(k)\La_{k-j},\cr}
$$
as claimed.
\cqfd

We now prove that 
$\underline \lambda \in {\cal Y}(\sharp\La , \La_{-u})$.
This amounts to showing that for each $k\ge0$,
\begin{equation}\label{flat_def}
p^\flat_k:=\Lambda_{k-u}+\sum_{j=0}^{n-1}\beta_j(k)\La_{k-j}
\end{equation}
is dominant. Assume the contrary. 
Then for some $i$, $p^\flat_{i}$ is dominant
and $p^\flat_{i-1}$ is not.
On using (\ref{ETADEF}), we have
$$
\mathalign{
p^\flat_{i-1}=p^\flat_i  -\epsilon_{i-1-u-\underline \lambda_i^\prime}
&=
\sum_{j=0}^{n-1} \beta_j(i)\Lambda_{i-j}+\Lambda_{i-u}
-\epsilon_{i-1-u-\underline \lambda_i^\prime}\cr
&=\sum_{j=0}^{n-1} \beta_j(i)\Lambda_{i-j}+\Lambda_{i-u}
-\Lambda_{i-u-\underline \lambda_i^\prime}+
\Lambda_{i-1-u-\underline \lambda_i^\prime}.\cr}
$$
Since $p^\flat_{i-1}$ is not dominant whereas $p^\flat_i$ is,
it follows that $\beta_{u+\underline \lambda_i^\prime}(i)=
-\delta^{(n)}_{0,\underline \lambda_i^\prime}$.
Comparing this expression with (\ref{flat_def}) at $k=i-1$
gives
$$
\sum_{j=0}^{n-1} \beta_j(i-1)\Lambda_{i-1-j}=
\sum_{j=0}^{n-1} \beta_j(i)\Lambda_{i-j}
+\Lambda_{i-u}-\Lambda_{i-1-u}
-\Lambda_{i-u-\underline \lambda_i^\prime}+
\Lambda_{i-1-u-\underline \lambda_i^\prime}\,,
$$
whereby, on applying the linear transformation $\La_j\mapsto \La_{i-1-j}$
for all $j$, we obtain:
$$
\mathalign{
p_{i-1}=
\sum_{j=0}^{n-1} \beta_j(i-1)\Lambda_j
&=\sum_{j=0}^{n-1} \beta_j(i)\Lambda_{j-1}
+\Lambda_{u-1}-\Lambda_{u}
-\Lambda_{u+\underline \lambda_i^\prime-1}+
\Lambda_{u+\underline \lambda_i^\prime}\cr
&=\sum_{j=0}^{n-1} \beta_j(i)\Lambda_{j-1}
-\epsilon_{u-1}+\epsilon_{u+\underline \lambda_i^\prime-1},\cr}
$$
so that
$$
p_i-p_{i-1}=
\sum_{j=0}^{n-1} \beta_j(i)\epsilon_{j-1}
+\epsilon_{u-1}-\epsilon_{u+\underline \lambda_i^\prime-1}.
$$
Since
$\beta_{u+\underline \lambda_i^\prime}(i)=
-\delta^{(n)}_{0,\underline \lambda_i^\prime}$,
the coefficient of $\epsilon_{u+\underline \lambda_i^\prime-1}$ on the right
side of this expression is negative, which contradicts
the fact that $p$ is a $\La$-path. 
Hence $\underline \lambda \in {\cal Y}(\sharp\La , \La_{-u})$.

Finally, it follows from Theorem \ref{diag_bij} that 
$\sharp\underline \lambda \in {\cal Y}(\La_u , \La)$
and that $\sharp\underline \lambda = \blambda$.
\cqfd
\begin{example}{\rm
To illustrate this theorem, let $n=4$ and $\La=\Lambda_1+2\Lambda_3$, and
consider the following highest-lift multipartition
$\blambda$ which has highest weight $\La$
$$
\smyoungd{\omit\phantom{\hrulefill}\cr
 \omit&\cr
 \omit\phantom{\hrulefill}\cr
 \omit&\cr
 \multispan{22}\hrulefill\cr
 &1&&2&&3&&0&&1&&2&&3&&0&&1&&&\omit&\cr
 \multispan{19}\hrulefill\cr
 &0&&1&&2&&3&&0&&1&&2&&3&\cr
 \multispan{17}\hrulefill\cr
 &3&&0&&1&&2&&3&&0&\cr
 \multispan{13}\hrulefill\cr
 &2&&3&&0&\cr
 \multispan{7}\hrulefill\cr
 &1&&2&&3&\cr
 \multispan{7}\hrulefill\cr
 &\cr
 \multispan{1}\hrulefill\cr
 }
\quad
\smyoungd{\multispan{22}\hrulefill\cr
 &3&&0&&1&&2&&3&&0&&1&&2&&3&&0&&\cr
 \multispan{21}\hrulefill\cr
 &2&&3&&0&&1&&2&&3&&0&&1&&2&&3&\cr
 \multispan{21}\hrulefill\cr
 &1&&2&&3&&0&&1&&2&\cr
 \multispan{13}\hrulefill\cr
 &0&&1&&2&&3&&0&&1&\cr
 \multispan{13}\hrulefill\cr
 &3&&0&&1&&2&&3&\cr
 \multispan{11}\hrulefill\cr
 &2&\cr
 \multispan{3}\hrulefill\cr
 &1&\cr
 \multispan{3}\hrulefill\cr
 &\cr
 \multispan{1}\hrulefill\cr
 }
\quad
\smyoungd{\multispan{22}\hrulefill\cr
 &3&&0&&1&&2&&3&&0&&1&&2&&&\omit&&\omit&\cr
 \multispan{17}\hrulefill\cr
 &2&&3&&0&&1&\cr
 \multispan{9}\hrulefill\cr
 &1&&2&&3&&0&\cr
 \multispan{9}\hrulefill\cr
 &0&&1&&2&\cr
 \multispan{7}\hrulefill\cr
 &\cr
 \multispan{1}\hrulefill\cr
 &\cr
 \multispan{1}\hrulefill\cr
 &\cr
 \multispan{1}\hrulefill\cr
 &\cr
 \multispan{1}\hrulefill\cr
 }
$$
Now construct the partition $\underline\lambda$ having, for each
$k=1,2,\ldots,$ the same number of rows of length $k$ as
does $\blambda$. On colouring $\underline \lambda$ with $-u=1$, we obtain
$$
\smyoungd{\multispan{22}\hrulefill\cr
 &1&&2&&3&&0&&1&&2&&3&&0&&1&&2&&\cr
 \multispan{21}\hrulefill\cr
 &0&&1&&2&&3&&0&&1&&2&&3&&0&&1&\cr
 \multispan{21}\hrulefill\cr
 &3&&0&&1&&2&&3&&0&&1&&2&&3&\cr
 \multispan{19}\hrulefill\cr
 &2&&3&&0&&1&&2&&3&&0&&1&\cr
 \multispan{17}\hrulefill\cr
 &1&&2&&3&&0&&1&&2&&3&&0&\cr
 \multispan{17}\hrulefill\cr
 &0&&1&&2&&3&&0&&1&\cr
 \multispan{13}\hrulefill\cr
 &3&&0&&1&&2&&3&&0&\cr
 \multispan{13}\hrulefill\cr
 &2&&3&&0&&1&&2&&3&\cr
 \multispan{13}\hrulefill\cr
 &1&&2&&3&&0&&1&\cr
 \multispan{11}\hrulefill\cr
 &0&&1&&2&&3&\cr
 \multispan{9}\hrulefill\cr
 &3&&0&&1&&2&\cr
 \multispan{9}\hrulefill\cr
 &2&&3&&0&\cr
 \multispan{7}\hrulefill\cr
 &1&&2&&3&\cr
 \multispan{7}\hrulefill\cr
 &0&&1&&2&\cr
 \multispan{7}\hrulefill\cr
 &3&\cr
 \multispan{3}\hrulefill\cr
 &2&\cr
 \multispan{3}\hrulefill\cr
 &\cr
 \multispan{1}\hrulefill\cr
 }
$$
We see that for each $k$, the colours at the ends of the
length $k$ rows of $\blambda$ are the $4$-complements of the colours
at the beginnings of the length $k$ rows of $\underline\lambda$.
Using Theorem \ref{mult_char}, we therefore conclude that
$\blambda\in{\cal Y}(\Lambda_3,\Lambda_1+2\Lambda_3)$.

If we colour $\underline\lambda$ instead with $-u\ne1$, we immediately
see that the colours at the ends of the
length $k$ rows of $\blambda$ are not the $4$-complements of the colours
at the beginnings of the length $k$ rows of $\underline\lambda$.
Therefore,
$\blambda\notin{\cal Y}(\Lambda_u,\Lambda_1+2\Lambda_3)$ if $u\ne3$.
\endex}
\end{example}

\def\AK{{\cal H}_m(\hbox{\boldmath$v$})}
\def\AKi{{\cal H}_m(\hbox{\bf i})}
\section{Ariki-Koike algebras}
\label{SECT3}
\subsection{Ariki-Koike algebras and affine Hecke algebras}
\label{SECT3.1}
The Ariki-Koike algebra ${\cal H}_m(v;u_0,\ldots,u_{l-1})$ is a deformation
of the group algebra of the complex reflection group
$G(l,1,m)=(\Z/l\Z)\wr \SG_m$, the wreath product of a cyclic
group of order $l$ with the symmetric group $\SG_m$. 
The group $G(l,1,m)$ can be realised as the group 
of monomial $m\times m$ matrices
whose entries are $l$th roots of unity. It is generated by the
permutation matrices of the elementary transpositions 
$$\sigma_i=(i,i+1), \qquad i=1,\ldots,m-1,$$
together with the matrix $\sigma_0 = {\rm diag\, }(\omega,1,\ldots,1)$, 
where $\omega$ is a primitive $l$th root of unity.
These generators satisfy
\begin{eqnarray}
&&\sigma_0\sigma_1\sigma_0\sigma_1=\sigma_1\sigma_0\sigma_1\sigma_0,\\
&&\sigma_0^l=1,\\ 
&&\sigma_0\sigma_i=\sigma_i\sigma_0, \qquad \mbox{   for $i\ge 2$,}
\end{eqnarray}
in addition to the usual Moore-Coxeter relations
for the generators of $\SG_m$.
\begin{definition} {\rm \cite{AK,BM}} The Ariki-Koike algebra
${\cal H}_m(v;u_0,u_1,\ldots,u_{l-1})$ is the unital associative $\C$-algebra
generated by $T_0,T_1,\ldots,T_{m-1}$ subject to the relations:
\begin{eqnarray}
&&T_iT_{i+1}T_i=T_{i+1}T_iT_{i+1},\quad 1\le i\le m-2,\label{EQ_T1}\\
&&T_iT_j=T_jT_i,\qquad \vert i-j\vert>1,\label{EQ_T2}\\
&&(T_i-v)(T_i+1)=0,\quad 1\le i\le m-1,\label{EQ_T3}\\
&&T_0T_1T_0T_1=T_1T_0T_1T_0,\\
&&(T_0-u_0)(T_0-u_1)\cdots(T_0-u_{l-1})=0,
\end{eqnarray}
where $v,u_0,u_1,\ldots,u_{l-1}$ are complex parameters.
\end{definition}

\noindent For convenience, we use $\AK$ to denote
${\cal H}_m(v;u_0,u_1,\ldots,u_{l-1})$.
\noindent In the case where $v=1$ and $u_j=\omega^j$, 
$\AK$ is isomorphic to the group algebra of $G(l,1,m)$.
For  $l=1$ and $l=2$, $\AK$ is
the Hecke algebra of types $A_{m-1}$ and $B_m$ respectively.

The group $G(l,1,m)$ is a quotient of the affine Weyl group
$\hat{W}_m=\Z\wr\SG_m$.  Similarly, $\AK$ is a quotient of
the affine Hecke algebra $\hat{H}_m(v)$, which is generated
by invertible elements
$T_1,\ldots, T_{m-1}$, $y_1,\ldots,y_m$ subject to 
the relations (\ref{EQ_T1}) (\ref{EQ_T2}) (\ref{EQ_T3}) 
of the Hecke algebra of type $A_{m-1}$ and to
\begin{eqnarray}
&&y_iy_j=y_jy_i, \qquad 1\le i,j\le m,\label{EQ_YY}\\
&&y_jT_i=T_iy_j,\qquad \mbox{for $j\not= i,i+1$},\label{EQ_TY2}\\
&&T_iy_iT_i=v\,y_{i+1}, \qquad 1\le i\le m-1 .\label{EQ_Y}
\end{eqnarray} 
Note that by (\ref{EQ_Y}), the generators $y_2,\ldots ,y_m$
are redundant, since
$$
y_{k}=v^{-k+1}T_{k-1}T_{k-2}\ldots T_1y_1T_1\ldots T_{k-2}T_{k-1},
\qquad 2\le k \le m.
$$
Also, it follows from Eqs.~(\ref{EQ_YY}) (\ref{EQ_Y}) that
\begin{equation}\label{EQ_TY1}
T_1y_1T_1y_1 = y_1T_1y_1T_1.
\end{equation}
In fact, $\hat{H}_m(v)$ may equally be defined as
the associative algebra generated by $T_1,\ldots ,T_{m-1}$ and $y_1$, 
subject to the relations (\ref{EQ_T1}) (\ref{EQ_T2}) (\ref{EQ_T3})
(\ref{EQ_TY1}) plus relation (\ref{EQ_TY2}) with $j=1$.
This shows that $\AK$ is the quotient of $\hat H_m(v)$ by the single relation
\begin{equation}
(y_1-u_0)(y_1-u_1)\cdots(y_1-u_{l-1})=0,
\end{equation} 
$T_0$ being taken equal to $y_1$.
The images of the $y_k$'s under the natural projection
of $\hat H_m(v)$ onto $\AK$ will still be denoted by $y_k$.
They generate a commutative subalgebra of $\AK$.
In the case where $l=1=u_0$, where $\AK$ reduces to the
Hecke algebra $H_m(v)$ of type $A_{m-1}$,
we have $y_1=T_0=1$, and the elements 
\begin{equation}
L_k = {y_k-1\over v-1}, \qquad 2\le k\le m
\end{equation}
are called the Jucys-Murphy elements of $H_m(v)$.

The affine Hecke algebra $\hat H_m(v)$ admits a faithful
representation by symmetrisation operators acting on
the ring of Laurent polynomials
$\C[x_1^{\pm 1}, \ldots, x_m^{\pm 1}]$ \cite{Lu,LS}. 
In this representation, the generators of $\hat H_m(v)$
act by
\begin{equation}
y_i f = x_i^{-1}f\,, 
\qquad T_i f = (v-1)\pi_i f +\sigma_i f\,,
\qquad f \in \C[x_1^{\pm 1}, \ldots, x_m^{\pm 1}]\,,
\end{equation}
where $\sigma_i$ is the transposition exchanging $x_i$ and
$x_{i+1}$ and $\pi_i$ is the isobaric divided difference
\begin{equation}
\pi_if = {x_i f - x_{i+1}\sigma_i(f)\over x_i -x_{i+1}} \,.
\end{equation}
Since symmetric expressions in the $x_i$ are scalars for these
operators, it is clear that the power-sums $y_1^k+\cdots+y_m^k$
are in the center of $\hat H_m(v)$ (they are in fact generators
of the center). In particular, 
\begin{equation}
c_m=y_1+\cdots+y_m
\end{equation}
acts as a scalar on every irreducible representation of $\hat H_m(v)$
or of its quotient $\AK$.

\subsection{Representation theory}
\label{SECT3.2}
The representation theory of $\AK$ has been studied  in \cite{AK,Ar1}.
We first have the following criterion of semi-simplicity.
\begin{theorem} {\rm\cite{Ar1}} \label{AKSStheorem}
The algebra $\AK$ is
semisimple if and only if the parameters satisfy
$$
v^d u_i \not = u_j \quad (\mbox{ for } i\not = j, d\in \Z,
|d|<n)
\quad \mbox{  and  }\quad [n]!\not = 0,  
$$
where
$[n]!=\prod_{j=1}^n[j]$ and 
$[j] = 1 + v +\cdots +v^{j-1}$.
\end{theorem}
In the semisimple case, the irreducible representations have
been constructed, and we have: 
\begin{theorem} {\rm\cite{AK}} \label{AKREPtheorem}
In the case where $\AK$ is semisimple,
the full set of non-equivalent irreducible representations
of $\AK$ is labelled by the set of $l$-tuples of partitions
$\blambda = (\lambda^{(0)},\ldots ,\lambda^{(l-1)})$ such
that $|\blambda|=\sum_j |\lambda^{(j)}|= m$.
\end{theorem}

\noindent
The irreducible $\AK$-module indexed by $\blambda$
will be denoted by $S(\blambda)$.
In \cite{AK}, a basis $\{b_\btau\}$ 
of  $S(\blambda)$ is constructed. It is labelled by standard multitableaux
of shape $\blambda$, that is, $l$-tuples $\btau=(\tau_0,\ldots,\tau_{l-1})$
of Young tableaux such that $\tau_i$ is of shape $\lambda^{(i)}$ and
each integer from $1$ to $m$  occurs exactly once.
The vectors $b_\btau$  are simultaneous eigenvectors
of the $y_k$, and there holds
\begin{equation}
y_k\, b_\btau = u_{s(k)}v^{c(k)}b_\btau\,,
\end{equation}
where $s(k)$ is the index of the component of $\btau$ in which
$k$ occurs, and $c(k)$ is the content of the box $k$ of $\lambda^{(s(k))}$.
It follows that $c_m$ acts on  $S(\blambda)$ as the scalar
\begin{equation}
c_\blambda
=
\sum_{j=0}^{l-1} u_{j} \left(\sum_{x\in \lambda^{(j)}} v^{c(x)}\right) \ .
\end{equation}

Suppose now that the parameters  are integral powers of $v$:
\begin{equation}\label{SpecPar}
u_j=v^{v_j} , \qquad 0\le j \le l-1 .
\end{equation}
Then, by Theorem~\ref{AKSStheorem}, $\H_m(\v)$ is not semisimple.
Let $M$ be any $\AK$ module on which $c_m$ acts as a scalar $c$,
and let $M\!\downarrow$ denote
the restriction of $M$ to ${\cal H}_{m-1}({\v})$.
Then $M\!\downarrow$ splits as the direct sum of the generalised 
eigenspaces of $c_{m-1}$. 
Since all $u_j$ are powers of $v$, the eigenvalues
of $c_{m-1}$ differ from those of $c_{m}$ by  just a power of $v$.
\begin{definition} {\rm \cite{Ar2}}
For $i\in\Z$, the $i$-restriction $M\!\downarrow_i$ of $M$
is defined as the generalised eigenspace of $c_{m-1}$ in
$M\!\downarrow$ corresponding to the eigenvalue $c-v^i$.
Similarly, the $i$-induced module $M\!\uparrow_i$ is the
$c_{m+1}$-generalised eigenspace of eigenvalue $c+v^i$ in the
induced module $M\!\uparrow$. 
\end{definition}
Let ${\cal G}_m(\v)=G_0(\AK)$ be the Grothendieck group of the
category of finite-dimensional  $\AK$-modules,
and 
$$
{\cal G}(\v)=\bigoplus_{m\ge 0} {\cal G}_m(\v) \ .
$$
As shown in \cite{LLT} for $l=1$, and in \cite{Ar2} for the general
case, one can define on ${\cal G}_\C(\v)=\C\otimes_\Z{\cal G}(\v)$
an action of the affine Kac-Moody algebra $\glchap_\infty$ if
$v$ is not a root of unity, and of $\slchap_n$ if $v$ is
a primitive $n$th root of unity, by setting
\begin{equation}
e_i [M] = [M\!\downarrow_i]\,,\qquad f_i[M] =[M\!\uparrow_i] \ .
\end{equation}
In both cases, ${\cal G}_\C(\v)$ is the level $l$ 
irreducible representation $V(\Lambda)$ with highest weight
$$
\Lambda=\Lambda_{v_0}+\cdots+\Lambda_{v_{l-1}} \ .
$$
Ariki's theorem can now be stated as follows:
\begin{theorem}{\rm \cite{Ar2}}\label{TH3-5}
Under the identification of ${\cal G}_\C(\v)$ with $V(\Lambda)$,
the basis of ${\cal G}_\C(\v)$ consisting of the classes
of the irreducible modules of the various $\AK$
is mapped to the canonical basis of $V(\Lambda)$.
\end{theorem}
Here, the canonical basis is understood in the sense of Lusztig.
It is the same as the global upper crystal basis of Kashiwara.
As a consequence, there is a one-to-one correspondence
between the simple $\AK$-modules and the vertices
of the crystal graph of $V(\La)$ whose associated weight
vectors have principal degree $m$.

Moreover, in the $l=1$ case and when $v$ is a primitive
root of unity of prime order $p$, it was shown in
\cite{LLT} that the crystal graph 
coincides with the $p$-good lattice of Kleshchev 
describing socles of restricted simple modules of
symmetric groups in characteristic $p$ \cite{Kl3}.

\subsection{The generalised Jantzen-Seitz problem}
\label{SECT3.3}
Henceforth, we assume that $v$ is a primitive $n$th root of
unity, and that the relations (\ref{SpecPar}) are satisfied.
Clearly, we may also assume that 
$0\le v_0\le \ldots \le v_{l-1}<n$.
We put ${\bf i}=(i_0,\ldots ,i_{n-1})$ where 
$i_k$ is the number of $v_j$ equal to $k$ and we write 
$\H_m({\bf i})$ for the Ariki-Koike algebra with this
choice of parameters.
We denote by $\La_{\bf i}=\sum_k i_k\La_k$
the corresponding $\slchap_n$-weight.

The Jantzen-Seitz problem consists in the determination of
those simple $\H_m({\bf i})$-mo\-dules which remain irreducible after
restriction to ${\cal H}_{m-1}({\bf i})$.
We shall translate this problem in the language of
crystal bases.
\begin{proposition}\label{REST}
Let $D$ be a simple $\AKi$-module, and let $b\in B(\La_{\bf i})$
be the corresponding vertex of the crystal graph of $V(\La_{\bf i})$.
Fix $k$ in $\{0,\ldots , l-1\}$.
The following three conditions on $j>0$ are equivalent:
\begin{quote}
{\rm (i)} $j$ is the smallest integer such that 
${D\!\downarrow\!_k}^{j+1}=0$.

{\rm (ii)} $\varepsilon_k(b) = j$.

{\rm (iii)} $[{D\!\downarrow\!_k}^{j}]=j!\,[D']$ for some simple
${\cal H}_{m-j}({\bf i})$-module $D'$.
\end{quote}
The equivalence between {\rm (i)} and {\rm (ii)} also holds for $j=0$.
\end{proposition}
\Proof 
Let $\{G(b) \mid b\in B(\La_{\bf i})\}$ denote the global 
upper crystal basis of $V(\La_{\bf i})$ (at $q=1$).
The fact that {\rm (i)} and {\rm (ii)} are equivalent 
follows from Lemma~5.1.1 (i), p.~470 of \cite{Ka3} 
and Theorem~\ref{TH3-5}.
Also, by Lemma~5.1.1 (ii) of \cite{Ka3}, if $\varepsilon_k(b) = j$
then 
$$
[{D\!\downarrow\!_k}^j]= e_k^{\,j}[D] = 
 e_k^{\,j}\,G(b) = j!\, G(\tilde e_k^{\,j} b)
= j!\, [D'],
$$
where $D'$ is simple by Theorem~\ref{TH3-5}.
Conversely assume {\rm (iii)}, \ie  
\begin{equation}\label{EQ*}
e_k^{(j)} G(b) = G(b'),
\end{equation}
where $e_k^{(j)} = e_k^{\,j}/j!$ is the $j$th divided power
of $e_k$ and $b'$ is the vertex corresponding to $D'$.
By repeated use of Eq. (5.3.8) of \cite{Ka3}, we see that
$$
e_k^{(j)} G(b) = {\varepsilon_k(b)\choose j} G(\tilde e_k^j b)
+ \ \mbox{other terms,}
$$
so that, if $j\not = 0$,  (\ref{EQ*}) implies that $j=\varepsilon_k(b)$.
\cqfd 

\begin{definition}
Let ${\bf j} = (j_0,\ldots ,j_{n-1}) \in \N^n$. 
We say that a simple $\AKi$-module $D$ satisfies the generalised
Jantzen-Seitz condition $\JS({\bf j})$ if and only if
$$
{D\!\downarrow\!_k}^{j_k+1}=0 \quad \mbox{\rm for $k=0,\ldots ,n-1$.}
$$
In this case, we write $D \in {\rm JS({\bf j})}$.
\end{definition}
It follows from Proposition~\ref{REST} that the set of irreducible
representations of $\AKi$ which restrict to irreducible
representations of ${\cal H}_{m-1}({\bf i})$ is
$$
\JS(1,0,\ldots ,0) \cup
\JS(0,1,0,\ldots ,0) \cup \cdots \cup \JS(0,\ldots ,0,1)\,.
$$
Proposition~\ref{REST} also implies that, in general, the
vertices $b$ of $B(\La_{\bf i})$ corresponding to 
simple $\AKi$-modules $D$ satisfying $\JS({\bf j})$ are
characterised by 
\begin{equation}
\varepsilon_k(b)\le j_k \quad \mbox{\rm for $k=0,\ldots ,n-1$.}
\end{equation}
For $\La\in P^+$, the number of those $b$ satisfying also
${\rm wt}(b)=\La - \La_{\bf j}$
is equal by Proposition~\ref{TPcorollary} to the tensor product
multiplicity $c_{\Lambda_{\bf j}\,\Lambda_{\bf i}}^\Lambda$
of $V(\Lambda)$ in
$V(\Lambda_{\bf j})\otimes V(\Lambda_{\bf i})$.
Therefore, the generating function of the number
of $\AKi$-modules $D$ satisfying $\JS({\bf j})$ is a sum of
branching functions of $\slchap_n$.

To state this result precisely, let us introduce some notation.
For $b\in B(\La_{\bf i})$, let $\deg b$ denote the homogeneous degree of
$b$, that is, $\deg b := -\langle \wt(b) \,,\,d  \rangle$.
If $b$ is labelled by a multipartition $\blambda$, then 
$\deg b = N^0(\blambda)$.
We write $b\in \JS({\bf j})$ if $G(b)=[D]$ 
for a module $D\in\JS({\bf j})$.
\begin{theorem}\label{TH3-8}
Let $l_{\bf i}$, $l_{\bf j}$ denote the levels of $\La_{\bf i}$
and $\La_{\bf j}$.
We have
$$
\sum_{b\in B(\La_{\bf i})\cap \JS({\bf j})}
z^{\deg b}
=
\sum_\La b^{\Lambda}_{\Lambda_{\bf j},\Lambda_{\bf i}}(z)\,,
$$
where $\La$ runs through the weights of $P^+_{l_{\bf i}+l_{\bf j}}$
congruent to $\La_{\bf i} + \La_{\bf j}$ modulo the $\Z$-lattice spanned
by the $\alpha'_i$.
\hfill $\Box$
\end{theorem}
\def\b{{\hbox{\boldmath $b$}}}
When one of the factors of the tensor product is of
level one, the branching functions $b^\Lambda_{\Lambda'\,\Lambda''}(z)$
can be explicitly evaluated in terms of theta functions
\cite{KP,JMO}, as we shall now recall. 

Let $\h=(\bigoplus_{i=0}^{n-1}\C\, h_i)\oplus\C d $
and $\bar\h = \bigoplus_{i=1}^{n-1}\C\, h_i$ denote
the Cartan subalgebras of $\slchap_n$ and $\Sl_n$, respectively. 
One identifies $\bar\h^*$ with the subspace
$\{\La\in \h^* \mid \La(c)=\La(d)=0\}$.
The linear map $\La \mapsto \bar\La$ defined by
$$
\bar\La_i=\Lambda_i-\Lambda_0,\qquad 
\bar\delta = 0,
$$
is a natural projection of $\h^*$ onto $\bar\h^*$.
Note that $\bar\alpha_i = \alpha'_i$.
The root lattice of $\Sl_n$ gets identified with
$\bar Q=\bigoplus_{i=1}^{n-1}\Z\bar\alpha_i$.
Let $(\cdot\,,\,\cdot)$ denote the usual bilinear form on $\h^*$
defined by
$$
(\La_i\,,\,\La_j) = \min (i,j) - {ij\over n}\,,\qquad
(\delta\,,\,\La_i) = 1,\qquad (\delta\,,\,\delta)=0.
$$
We write $|\La|^2=(\La\,,\,\La)$.
Recall that the Weyl group of $\Sl_n$ is isomorphic to $\SG_n$,
the elementary transposition $\sigma_i$ being represented
by the orthogonal reflection of $\bar\h^*$ with respect
to $\bar\alpha_i$.

For $\mu\in\bar\h$ and $m\in\Z_{>0}$, define the {\em theta function}
\begin{equation}
\Theta_{\mu,m}(z)=\sum_{\alpha\in\bar Q}
z^{ {m\over 2} \left| \alpha+ {1\over m}\mu \right|^2} \ ,
\end{equation}
and let $\eta(z)=z^{1/24}\prod_{k\ge 1}(1-z^k)$ be the Dedekind
modular function. Then, for a fundamental weight $\LPP$, 
and for $\La \equiv \LP+\LPP \mod \bar Q$,
Eq.~(4.17) of \cite{JMO} reads
\begin{equation}\label{eq:bth}
\b^\Lambda_{\Lambda'\,\Lambda''}(z)
=
\eta(z)^{-(n-1)}\sum_{w\in\SG_n} {\rm sgn}(w)
\Theta_{ -L(\bar\Lambda'+\bar\rho)+(L-1)w(\bar\Lambda+\bar\rho)\,
,\,L(L-1)}(z) \ ,
\end{equation}
where $L=n+(\delta,\lambda)$ and $\rho = \sum_{i=0}^{n-1}\La_i$. 
The branching function $b^\Lambda_{\Lambda'\,\Lambda''}(z)$ of
Eq.~(\ref{BRFU}) differs from 
$\b^\Lambda_{\Lambda'\,\Lambda''}(z)$ by a power of $z$,
namely  
$$
\b^\Lambda_{\Lambda'\,\Lambda''}(z)
=z^{\gamma(\Lambda',\Lambda'';\Lambda)}b^\Lambda_{\Lambda'\,\Lambda''}(z)
$$
where 
\begin{equation}\label{eq:gamma}
\gamma(\Lambda',\Lambda'';\Lambda)
=
{|\Lambda'+\rho|^2 \over 2(l'+n)} +
{|\Lambda''+\rho|^2 \over 2(l''+n)} -
{|\Lambda+\rho|^2 \over 2(l+n)} -
{|\rho|^2 \over 2n} 
\end{equation}
and $l, l', l''$ denote the respective levels of 
$\La, \La', \La''$.
\begin{example}{\rm
Take $n=2$, $\Lambda'=\Lambda_0+\Lambda_1$ and $\Lambda''=\Lambda_0$.
Then $L=5$, and the nonzero branching functions correspond to 
$\Lambda=2\Lambda_0+\Lambda_1$ and $\Lambda=3\Lambda_1$.
Equation (\ref{eq:bth}) gives
$$
\b^{2\Lambda_0+\Lambda_1}_{\Lambda_0+\Lambda_1,\Lambda_0}(z)=
\eta(z)^{-1}\left[
\Theta_{-2\bar\Lambda_1,20}(z)-\Theta_{-18\bar\Lambda_1,20}(z)\right] \ .
$$
Here,
$$
\Theta_{-2\bar\Lambda_1,20}(z)=\sum_{k\in\Z}
z^{10\left|k\alpha_1-{1\over 20}\alpha_1\right|^2}
=\sum_{k\in\Z}
z^{20\left(k-{1\over 20}\right)^2} \ ,
$$
$$
\Theta_{-18\bar\Lambda_1,20}(z)=\sum_{k\in\Z}
z^{10\left|k\alpha_1-{9\over 20}\alpha_1\right|^2}
=\sum_{k\in\Z}
z^{20\left(k-{9\over 20}\right)^2}  \ ,
$$
so that
$$
\b^{2\Lambda_0+\Lambda_1}_{\Lambda_0+\Lambda_1,\Lambda_0}(z)=
z^{1/20}\eta(z)^{-1}\sum_{k\in\Z}
\left(
z^{20k^2-2k} - z^{20k^2-18k+4}
\right)
$$
$$
= z^{1/120}(1+z+2z^2+3z^3+4z^4+6z^5 +8z^6+11z^7+ \cdots\, ) \ ,
$$
and since (\ref{eq:gamma}) gives $\gamma={1\over 120}$, we have
finally
$$
b^{2\Lambda_0+\Lambda_1}_{\Lambda_0+\Lambda_1,\Lambda_0}(z)=
1+z+2z^2+3z^3+4z^4+6z^5 +8z^6+11z^7+ \cdots\, \ .
$$
Likewise,
$$
\b^{3\Lambda_1}_{\Lambda_0+\Lambda_1,\Lambda_0}(z)=
\eta(z)^{-1}\left[
\Theta_{6\bar\Lambda_1,20}(z)-\Theta_{-26\bar\Lambda_1,20}(z)\right]
$$
$$
=z^{9/20}\eta(z)^{-1}
\sum_{k\in\Z}\left(
z^{20k^2+6k}-z^{20k^2-26k+8}\right)
$$
$$
=z^{49/120}(1+z+z^2+2z^3+3z^4+ 4z^5+6z^6+8z^7+\cdots\,) \ .
$$
Here, (\ref{eq:gamma}) yields $\gamma=-{71\over 120}$
so that
$$
b^{3\Lambda_1}_{\Lambda_0+\Lambda_1,\Lambda_0}(z)=
z+z^2+z^3+2z^4+3z^5+4z^6+6z^7+\cdots\ .
$$

In terms of Ariki-Koike algebras, using Theorem~\ref{TH3-8},
we have found that
the generating function for the number of $\H_m(1,0)$-modules
satisfying ${\rm JS}(1,1)$, 
(or equivalently for the number of ${\cal H}_m(1,1)$-modules
satisfying ${\rm JS}(1,0)$) is equal to
\begin{eqnarray*}
&&b^{2\Lambda_0+\Lambda_1}_{\Lambda_0+\Lambda_1,\Lambda_0}(z)
+
b^{3\Lambda_1}_{\Lambda_0+\Lambda_1,\Lambda_0}(z) \\
&&\qquad\qquad\quad = \
\prod_{i>0}{1\over 1-z^i}
\sum_{k\in\Z}
(z^{20k^2-2k}+z^{20k^2+6k+1}-z^{20k^2-18k+4}
-z^{20k^2-26k+9}) \\
&&\qquad\qquad\quad = \ 
1+2z+3z^2+4z^3+6z^4+9z^5+12z^6+17z^7\cdots\ .
\end{eqnarray*}
\endex}
\end{example}
%
\subsection{Labelling of modules}
\label{SECT3.4}
Let $\La_{\bf i} \in P_l^+$.
As explained in Section~\ref{SECT2.1}, there are (at least)
three ways of labelling the vertices of the crystal graph
of $V(\La_{\bf i})$.
One may use paths $p\in \P(\La_{\bf i})$, their highest-lifts
$\blambda \in {\cal Y}(\La_{\bf i})$, or the set
${\cal M}(\La_{\bf i})$ of multipartitions coming from the
second $q$-deformation $\F'_q(\La_{\bf i})$ 
of the level $l$ Fock space $\F(\La_{\bf i})$
(see Note~\ref{NoteFockSpace}).

From the point of view of Ariki-Koike algebras, there
is a canonical way of labelling the irreducible
representations of $\AKi$. This is as follows.
Under the specialisation of parameters of Section~\ref{SECT3.3}, 
the modules
$S(\blambda)$ are in general reducible.
An important property of $S(\blambda)$ is that it
is endowed with a natural bilinear form
compatible with the action of $\AKi$.
Let $\hbox{rad}\,S(\blambda)$ denote the radical of
this form. It was proved by Graham and Lehrer that
$S(\blambda)/\,\hbox{rad}\,S(\blambda)$ is either
a simple module or 0, and that all simple modules
arise this way \cite{GL}. (This generalises
a similar result of Dipper and James for type $A$ Hecke algebras
\cite{DJ1} and Dipper, James and Murphy  for type $B$ \cite{DJM}.)
Then Mathas, building on Ariki's theorem, characterised
the set of multipartitions $\blambda$ such that
$S(\blambda)/\,\hbox{rad}\,S(\blambda) \not = 0$
and proved that it coincides with ${\cal M}(\La_{\bf i})$
\cite{Ma}.
(Strictly speaking, Mathas' labels are obtained by conjugating and
switching the components of the $\blambda \in {\cal M}(\La_{\bf i})$.
This arises from a different convention for labelling the
$S(\blambda)$. We shall ignore here such minor modifications.)
This makes it natural to define 
$$
D(\blambda) = S(\blambda)/\,\hbox{rad}\,S(\blambda),\qquad
(\blambda \in {\cal M}(\La_{\bf i})).
$$

Unfortunately, it appears difficult to give a non-recursive description
of ${\cal M}(\La_{\bf i})$, whereas Proposition~\ref{TREVOR}
provides such a simple characterisation of ${\cal Y}(\La_{\bf i})$.
This prompts us to use ${\cal Y}(\La_{\bf i})$ as an alternative
set of labels for irreducible representations of $\AKi$.
The bijection between the two labellings may be obtained 
by following a sequence of arrows
back to the highest weight vertex in the crystal graph labelled
by ${\cal Y}(\La_{\bf i})$, and
then applying the reversed sequence to the highest weight vertex
of the crystal graph labelled by ${\cal M}(\La_{\bf i})$.
We denote the image of $\blambda\in{\cal Y}(\Lambda_{\bf i})$
under this bijection by $\tilde\blambda\in{\cal M}(\Lambda_{\bf i})$.
Then, for $\blambda\in{\cal Y}(\Lambda_{\bf i})$, define
$\tilde D(\blambda)=D(\tilde\blambda)$.
\begin{figure}[t]
\begin{center}
\leavevmode
\epsfxsize =13.5cm
\epsffile{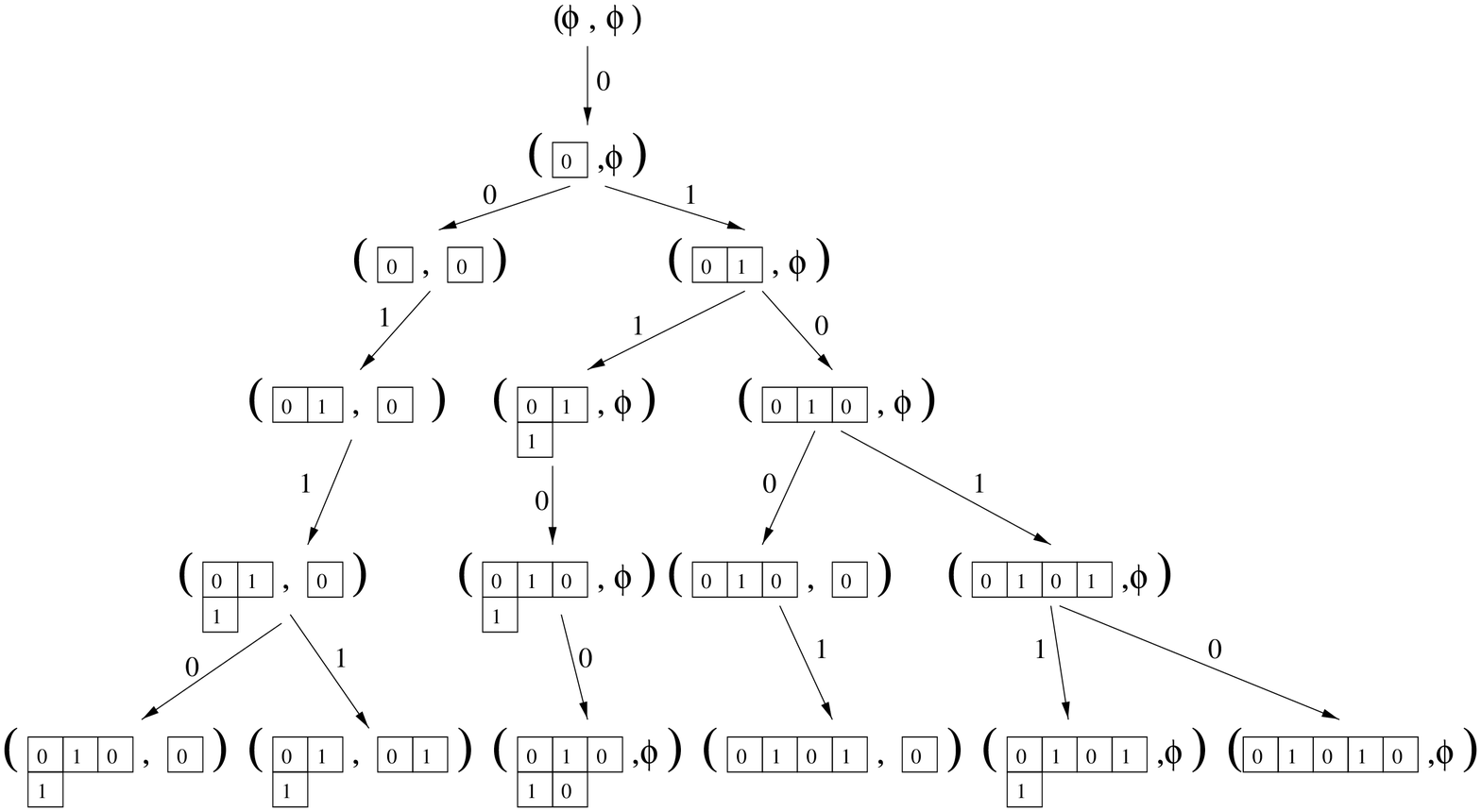}
\end{center}
\caption{\label{FIG2} The crystal graph the $U_q(\slchap_2)$-module 
$V_q(2\La_0)$ labelled by ${\cal M}(2\La_0)$}
\end{figure}
\begin{example}\rm
The crystal graph of the $U_q(\slchap_2)$-module $V_q(2\La_0)$
labelled by ${\cal M}(2\La_0)$ is shown in Figure~\ref{FIG2}.
By comparing with Figure~\ref{FIG1}, we get the following relations
between the two labellings of 
$\H_m(2,0)$-modules:
$$
\tilde D((2),(2)) = D((2,1),(1)), \qquad 
\tilde D((3),(2)) = D((3,1),(1)).
$$
The other modules for $m\le 5$ have the same label in both
systems of indexation.
\endex
\end{example}
We can now apply the results of Section~\ref{SECT2.3} and
Section~\ref{SECT2.4} and state:
\begin{corollary}\label{CORO}
$\strut$
\begin{quote}
{\rm (i)} The $\AKi$-module $\tilde D(\blambda)$ satisfies
the condition $\JS({\bf j})$ if and only if 
$\blambda \in {\cal Y}(\La_{\bf j},\La_{\bf i})$.

{\rm (ii)} The $\AKi$-module $\tilde D(\blambda)$ restricts
to a simple ${\cal H}_{m-1}({\bf i})$-module
if and only if $\blambda$ satisfies the conditions
of Theorem~{\rm\ref{mult_char}} for some integer $u$ between
$0$ and $n-1$.

{\rm (iii)} The map $\tilde D(\blambda) \mapsto \tilde D(\sharp\blambda)$ 
establishes a one-to-one correspondence between the simple $\AKi$-modules
satisfying $\JS({\bf j})$ and the simple 
${\cal H}_m(\sharp{\bf j})$-modules satisfying $\JS(\sharp{\bf i})$,
where $\sharp {\bf i}=(i_0,i_{n-1},\ldots ,i_1)$. 
\end{quote}
\end{corollary}
\begin{example}\rm
Let us determine the simple $\H_5(2,0)$-modules $D \in \JS(1,0) \cup \JS(0,1)$.
Recall that the Hecke algebra $H_m$ of type $B_m$ depends on two parameters
$q$ and $Q$ in the notation of \cite{DJM}, and coincides with the
AK-algebra $\H_m(v;u_0,u_1)$ where $v=q$, $u_0=-1$ and $u_1=Q$.
Therefore, if $q=Q=-1$, $H_m$ coincides with $\H_m(0,2)$, which
is isomorphic to $\H_m(2,0)$ by sending $T_0$ to $-T_0$.
Thus, this question is equivalent to 
finding the irreducible representations of the 
Hecke algebra of type $B_5$ with both parameters $q$ and $Q$ equal to $-1$,
which restrict to irreducible representations of the Hecke algebra 
of type $B_4$ (with the same parameters).

By Proposition~\ref{TREVOR}, or by inspecting the graph of Figure~\ref{FIG1},
one determines the following list of simple $\H_5(2,0)$-modules:
$$
\tilde D((3),(2)),\ \
\tilde D((2,1),(2)),\ \
\tilde D((3,2),\emptyset),\ \
\tilde D((4),(1)),\ \
\tilde D((4,1),\emptyset),\ \
\tilde D((5),\emptyset).
$$
Among them, by Corollary~\ref{CORO}~(i), the modules 
$\tilde D((3),(2))$ and $\tilde D((5),\emptyset)$ satisfy
$\JS(1,0)$ and the module $\tilde D((4),(1))$ satisfies
$\JS(0,1)$. This can also be checked using Figure~\ref{FIG1}.

Moreover, by Corollary~\ref{CORO}~(iii),
$\tilde D((3),(2))$ and $\tilde D((5),\emptyset)$
are in correspondence with the $\H_5(1,0)$-modules 
$\tilde D(3,2)$ and $\tilde D(5)$, respectively, and
$\tilde D((4),(1))$ corresponds to the $\H_5(0,1)$-module 
$\tilde D(4,1)$. These are the modules of each algebra satisfying 
the condition $\JS(2,0)$. 
Note that $\H_5(1,0)$ and $\H_5(0,1)$ are both isomorphic
to the Hecke algebra of type $A_4$ at $q=-1$.
\endex
\end{example}

\begin{appendix}
	
\section{Proofs of Propositions~\ref{TREVOR} and \ref{EndColLemma}
and Lemma \ref{WtLemma}}

\noindent
{\it Proof of} \ref{TREVOR}.
We shall prove that \ref{TREVOR} is equivalent to Proposition~3.4
of \cite{JMMO}.
Let $\blambda=(\lambda^{(0)},\ldots,\lambda^{(l-1)})$.
We put
\begin{equation} \label{TIJDEF}
t_{jk}=v_j-\lambda^{(j)\prime}_{k+1},\qquad (k\ge0,\ 0\le j<l),
\end{equation}
and we extend this definition to all $j\in\Z$ by setting
\begin{equation}\label{period}
t_{j+l,k}=t_{j,k}+n.
\end{equation}
Then $\blambda$ is a cylindrical multipartition of highest weight
$\La$ if and only if
\begin{equation}\label{CYL}
t_{j,k}\le t_{j+1,k}, \qquad (k\ge 0,\ j\in \Z).
\end{equation}
By  Proposition~3.4 of \cite{JMMO}, $\blambda \in {\cal Y}(\La)$
if and only if for each $k\ge0$, there exists an integer $j_*=j_*(k)$ such
that
\begin{equation}\label{HL}
t_{j_*+1,k}>t_{j_*,k+1}.
\end{equation}
Let us show that (\ref{CYL}) and (\ref{HL}) imply the characterisation 
given in Proposition \ref{TREVOR}.
For $j_*<j\le j_*+l$,
$$
t_{j,k}\ge t_{j_*+1,k}> t_{j_*,k+1}= t_{j_*+l,k+1}-n\ge t_{j,k+1}-n,
$$
whence $t_{j,k}>t_{j,k+1}-n$ holds for all $j$ by (\ref{period}).
Thus $\lambda^{(j)\prime}_{k+1}>\lambda^{(j)\prime}_k-n$, so that
each partition $\lambda^{(j)}$ is $n$-regular.
Therefore, $\lambda^{(j)}$  has $t_{j,k}-t_{j,k-1}<n$ rows of length $k$.
Using (\ref{TIJDEF}), we find that the colours at the end of these
rows are the ${}\mod n$ residues of
\begin{equation}\label{end_colours}
t_{j,k-1}+k, t_{j,k-1}+k+1, t_{j,k-1}+k+2, \ldots, t_{j,k}+k-1,
\end{equation}
whereupon those at the beginning of these rows
are the ${}\mod n$ residues of 
\begin{equation}\label{beg_colours}
t_{j,k-1}+1, t_{j,k-1}+2, t_{j,k-1}+3, \ldots, t_{j,k}.
\end{equation}
It follows that the
colours appearing at the beginnings of all the rows of length $k$
in $\blambda$
are the $\mod n$ residues of (\ref{beg_colours}) as $j$ varies
over $j_*(k-1)<j\le j_*(k-1)+l$.  Set $j_*=j_*(k-1)$.
For $0\le r<n$, let $r^\prime\in\Z$ be the unique value such that
$r^\prime\bmod n=r$ and such that $t_{j_*,k}<r^\prime\le t_{j_*+l,k}$.
Since $t_{j_*+l,k}=t_{j_*,k}+n<t_{j_*+1,k-1}+n$, this $r^\prime$ is
the only possible value appearing in (\ref{beg_colours})
for $j_*<j\le j_*+l$ that gives rise to $r$.
Let $\theta_{r^\prime,k}$ be the smallest value of $\theta$ such
that $t_{\theta,k}\ge r^\prime$. Then the number of times that
$r^\prime$ appears in (\ref{beg_colours}) for $j_*<j\le j_*+l$,
is equal to 
\begin{equation}\label{number_ip}
\psi_{r^\prime+k-1}=\theta_{r^\prime,k-1}-\theta_{r^\prime,k},
\end{equation}
for all $r^\prime\in\Z$.
(Note that $\theta_{r^\prime,k-1}\ge\theta_{r^\prime,k}$ for all
$r^\prime$.)
With this definition, the property (\ref{period}) of $t_{jk}$ implies
that this number is also given by $\psi_{r+k-1}=\psi_{r^\prime+k-1}$.

In the case $r^\prime=t_{j_*+1,k-1}$, note that
$\theta_{r^\prime,k-1}=\theta_{r^\prime,k}$ so that 
$\psi_{r^\prime+k-1}=0$ implying that $t_{j_*+1,k-1}\bmod n$ does not
appear at the beginning of a row of length $k$.
Furthermore, on comparing (\ref{end_colours}) and (\ref{beg_colours}),
we find that
the number of $r$-nodes at the
end of length $k$ rows of $\blambda$ is given by $\psi_r$.
In particular, for $r^\prime=t_{j_*+1,k-1}$, the colour 
$(r^\prime+k-1)\,\mod n$ does not appear at the end of a row of length $k$.

Hence (\ref{CYL}) and (\ref{HL}) imply the characterisation
given in Proposition \ref{TREVOR}.
This reasoning may be reversed, thus showing the converse. \cqfd

\noindent
{\it Proof of} \ref{EndColLemma}.1.
Let $p\in\P(\La)$ and $\blambda = \blambda(p)$.
Let $\psi_r(k)$ be the number of $r$-nodes at the right ends
of length $k$ rows of $\blambda$.
{}From (\ref{ETADEF}) and (\ref{TIJDEF}), we obtain
\begin{equation}\label{ETATIJ}
\eta_k=\sum_{j=0}^{l-1} \epsilon_{t_{j,k}+k},
\end{equation}
so that
\begin{eqnarray*}
\eta_{k-1}-\sigma(\eta_k)
&=&
\sum_{j=0}^{l-1} \left(
  \epsilon_{t_{j,k-1}+k-1} - \epsilon_{t_{j,k}+k-1} \right)\\
&=&
\sum_{j=0}^{l-1} \left(
  \alpha_{t_{j,k-1}+k}^\prime
  + \alpha_{t_{j,k-1}+k-1}^\prime
  + \cdots
  + \alpha_{t_{j,k}+k-1}^\prime \right).
\end{eqnarray*}
Since the subscripts here are the values in (\ref{end_colours}),
we immediately obtain
\begin{equation}\label{end_col_res}
\eta_{k-1}-\sigma(\eta_k)
=\sum_{r=0}^{n-1} \psi_r(k) \alpha_r^\prime.
\end{equation}
\cqfd

\noindent
{\it Proof of} \ref{EndColLemma}.2.
For $K\ge0$, define a path $p^{(K)}$
by $p^{(K)}_K=p^{\phantom{(}}_K$ and 
\begin{equation}\label{K_path}
p^{(K)}_{k+1}-p^{(K)}_k
=\sum_{j=0}^{l-1} \epsilon_{t_{j,K}+k}, \qquad (k\ge 0).
\end{equation}
Then $p^{(K)}$ is the ground state of
${\cal P}(\La^{(K)})$ where
$$
\La^{(K)}=\sum_{j=0}^{l-1} \Lambda_{t_{j,K}}\,.
$$
This also implies that $p^{(K)}_{K+1}=p^{\phantom{(}}_{K+1}$, so that
$p^{(K)}$ has its $K$th and $(K+1)$th points in common with $p$.
For $k,K>0$, it follows from (\ref{K_path}) that,
\begin{equation}\label{K_diff_ind}
p^{(K-1)}_{k-1}-p^{(K)}_{k-1}
=p^{(K-1)}_{k}-p^{(K)}_{k}
+\sum_{j=0}^{l-1} \epsilon_{t_{j,K}+k-1}
-\sum_{j=0}^{l-1} \epsilon_{t_{j,K-1}+k-1}.
\end{equation}
Together, (\ref{ETATIJ}) and (\ref{end_col_res}) imply that
\begin{equation}\label{K_diff}
\sum_{j=0}^{l-1} \epsilon_{t_{j,K}+K-1}
-\sum_{j=0}^{l-1} \epsilon_{t_{j,K-1}+K-1}
=\sigma(\eta_K) -\eta_{K-1}
=-\sum_{r=0}^{n-1} \psi_r(K) \alpha^\prime_r,
\end{equation}
where $\psi_r(K)$ is the number of $r$-nodes
appearing at the end of a length $K$ row of $\blambda$.
Shifting the subscripts in (\ref{K_diff}) and combining with
(\ref{K_diff_ind}) yields:
\begin{equation}\label{interp_step}
p^{(K-1)}_{k-1}-p^{(K)}_{k-1}
=p^{(K-1)}_{k}-p^{(K)}_{k}
-\sum_{r=0}^{n-1} \psi_r(K) \alpha^\prime_{r-K+k},\qquad (k,K>0).
\end{equation}
Then, for $k\le K$, repeated use of (\ref{interp_step}) yields:
\begin{eqnarray*}
p^{(K-1)}_{k-1}-p^{(K)}_{k-1}
&=&
p^{(K-1)}_{K}-p^{(K)}_{K}
  -\sum_{j=k}^K\sum_{r=0}^{n-1} \psi_r(K) \alpha_{r-K+j}^\prime\\
&=&
-\sum_{r=0}^{n-1} m_r(K,k) \alpha^\prime_r,
\end{eqnarray*}
(using $p^{(K-1)}_K=p^{\phantom{(}}_K=p^{(K)}_K$)
where $m_r(K,k)$ is the multiplicity of the colour $r$
in the length $K$ rows of $\blambda$ in or to the
right of the $k$th column.
Hence, if $K_*$ is the length of $p$ so that $p^{(K_*)}$
is the ground state of ${\cal P}(\La)$, then for $k\le K_*$,
$$
p^{\phantom{(}}_{k-1}-p^{(K_*)}_{k-1}=
p^{(k-1)}_{k-1}-p^{(K_*)}_{k-1}=
-\sum_{i=k}^{K_*}\sum_{r=0}^{n-1} m_r(i,k) \alpha^\prime_r=
-\sum_{r=0}^{n-1} m_r(k) \alpha^\prime_r,
$$
where $m_r(k)$ is the number of $r$-nodes
in $\blambda$ in or to the right of the $k$th column.
For $k>K_*$, both sides of this expression are clearly 0.
\cqfd

\noindent
{\it Proof of} \ref{WtLemma}.
Let $p\in\P(\La)$ and $\blambda = \blambda(p)$.
Proposition \ref{EndColLemma}(2) implies that
$$
p_0 = \Lambda + N^0(\blambda)\delta
    -\sum_{r=0}^{n-1} N^r(\blambda)\alpha_r.
$$
Once it is established that $N^0(\blambda)=E(p)$, we immediately
obtain the desired result from the definitions
(\ref{MultiWtDef}) and (\ref{PathWtDef}):
$$
\wt(p) = p_0-E(p)
=\Lambda-\sum_{r=0}^{n-1} N^r(\blambda)\alpha_r
=\wt(\blambda).
$$
The fact that $N^0(\blambda)=E(p)$ is not immediate.
The details of its verification may be extracted from
Prop. 5.6 of \cite{DJKMO}.
\cqfd

\def\Mu{{\rm M}}

\section{Tensor products and the $\sharp$ involution}

The $\sharp$ involution allowed us to map
the highest weight vertices of the
crystal graph of $V_q(\LP)\otimes V_q(\LPP)$ to those
of $V_q(\sharp\LPP)\otimes V_q(\sharp\LP)$. In this appendix,
we sketch the construction of a similar map at the level of 
the $U_q(\slchap_n)$-modules.

The symmetry $i\leftrightarrow n-i$ of the Dynkin diagram
of type $A_{n-1}^{(1)}$ induces an involutive automorphism
of $U_q(\slchap_n)$
given by
$$
e_i^\sharp=e_{-i}\,,\qquad 
f_i^\sharp=f_{-i}\,,\qquad
(q^{h_i})^\sharp=q^{h_{-i}}\,,\qquad
(q^d)^\sharp=q^d \,.
$$
Given a multipartition $\blambda=(\lambda^{(0)},\ldots,\lambda^{(l-1)})$
we define
$
\blambda'=({\lambda^{(l-1)}}',\ldots,{\lambda^{(0)}}')
$
where ${}'$ is the conjugation of partitions.
If $\blambda$ labels a vector $v_\blambda$ of $\F_q(\Lambda)$,
so that the nodes of the main diagonal of $\lambda^{(j)}$
are coloured with $v_j \mod n$, then we regard
$\blambda'$ as labelling a vector $v_{\blambda'}$ of 
$\F_q(\sharp\Lambda)$ and we colour the nodes of the main 
diagonal of $\lambda'^{(j)}$ with $-v_{l-1-j} \mod n$.
The following lemma is a generalisation to level $l>1$ of
Lemma~7.5 of \cite{LLT}.
\begin{lemma}\label{le:B}
Let $x\mapsto x'$ be the semi-linear map 
from $\F_q(\Lambda)$ to $\F_q(\sharp\Lambda)$ defined by
$$
\left(\sum_{\sblambda} \phi_\sblambda(q)\, v_\sblambda\right)'
=\sum_{\sblambda} \phi_\sblambda(q^{-1})\,v_{\sblambda'}\ .
$$
Then, one has
$$
f_i^\sharp (u') = (q^{-1-h_i} f_iu)',\quad
e_i^\sharp (u')=(q^{h_i-1}e_i u)', \qquad (u\in\F_q(\La),\ 0\le i\le n-1).
$$
\end{lemma}
\Proof
Let $\<\cdot\,,\,\cdot\>$ be the scalar product on $\F_q(\Lambda)$ for which
the natural basis $\{v_\blambda\}$ is orthonormal. 
Using the formulae of Theorem~\ref{TH2-5}, one can
check that, for $u,v \in \F_q(\La)$,
\begin{eqnarray*}
\<q^hu\,,\,v\> &=& \<u\,,\,q^hv\>,\qquad\quad (h\in\h),\\
\<f_iu\,,\,v\> &=& \<u\,,\,q^{h_i-1}e_iv\>, \qquad (0\le i\le n-1),\\
\<e_iu\,,\,v\> &=& \<u\,,\,q^{-1-h_i}f_iv\>, \qquad (0\le i\le n-1) .
\end{eqnarray*}
Let us prove that 
$$
f_i^\sharp (u') = (q^{-1-h_i} f_iu)',
\qquad (u\in\F_q(\La),\ 0\le i\le n-1).
$$
By linearity, it is sufficient to show that
\begin{equation}
\<f_i^\sharp v_{\sblambda'}\,,\, v_{\sbmu'}\> =
\< (q^{-1-h_i}f_iv_{\sblambda})'\,,\,v_{\sbmu'}\>
\end{equation}
for all basis vectors $v_\sblambda$ of $\F_q(\La)$ and $v_{\sbmu'}$ 
of $\F_q(\sharp\La)$.
Using again Theorem~\ref{TH2-5} for 
$\F_q(\Lambda)$ and $\F_q(\sharp\Lambda)$, we see that
$$
\<f_i^\sharp v_{\sblambda'}\,,\,v_{\sbmu'}\>
=
\overline{\<v_\sblambda\,,\,e_iv_\sbmu\>}\,,
$$
where, for $\phi(q)\in\Q(q)$, we set $\overline{\phi(q)}=\phi(q^{-1})$.
Then,
$$
\overline{\<v_\sblambda\,,\,e_iv_\sbmu\>}=
\overline{ \<q^{-1-h_i}f_iv_\sblambda\,,\,v_\sbmu\>}
=
\<(q^{-1-h_i}f_iv_\sblambda)'\,,\,v_{\sbmu'}\>
$$
whence the result. 
The formula for $e_i^\sharp (u')$ is proved similarly. \cqfd

Consider now the tensor product $\F_q(\LP)\otimes \F_q(\LPP)$
with the $U_q(\slchap_n)$-module structure determined by
the comultiplication (\ref{COM}).
Define a semi-linear map ${}'$ from
$\F_q(\LP)\otimes \F_q(\LPP)$ to $\F_q(\sharp\LPP)\otimes \F_q(\sharp\LP)$
by
\begin{equation}
u\otimes v\mapsto (u\otimes v)'=v'\otimes u' \ .
\end{equation}
Introducing  on
$\F_q(\LP)\otimes \F_q(\LPP)$ the scalar product $\<\cdot\,,\,\cdot\>$ 
for which the tensors $v_\blambda\otimes v_\bmu$ form an orthonormal basis,
and arguing as in the proof of Lemma \ref{le:B}, we obtain:
\begin{lemma}\label{B2}
Let $w\in\F_q(\LP)\otimes \F_q(\LPP)$. Then,
$$
f_i^\sharp(w')=(q^{-1-h_i}f_iw)'\,,\quad\mbox{and}\ \
e_i^\sharp(w')=(q^{h_i-1}e_iw)' \ .
$$
In particular, if $w$ is a highest weight vector of 
weight $\La$ in $\F_q(\LP)\otimes \F_q(\LPP)$, 
then $w'$ is a highest weight vector of weight $\sharp\La$ in 
$\F_q(\sharp\LPP)\otimes \F_q(\sharp\LP)$.
\end{lemma}
Let $V_q(\Lambda)\subset\F_q(\Lambda)$ be the irreducible component
generated by the vacuum vector $\emptyset_\Lambda$ of $\F_q(\Lambda)$.
The tensor product $V_q(\LP)\otimes V_q(\LPP)$ is the submodule of
$\F_q(\LP)\otimes \F_q(\LPP)$ generated by the action of
$U_q(\slchap_n)\otimes U_q(\slchap_n)$ on 
$\emptyset_{\LP}\otimes\emptyset_{\LPP}$.
From Lemma~\ref{le:B}, we see that $V_q(\LP)\otimes V_q(\LPP)$ is
mapped by ${}'$ onto $V_q(\sharp\LPP)\otimes V_q(\sharp\LP)$,
and from Lemma~\ref{B2}, we deduce 
that the map ${}'$ respects the decomposition
into $U_q(\slchap_n)$-modules. 

In particular, highest weight
vectors of weight $\La$ of 
$V_q(\LP)\otimes V_q(\LPP)$ are mapped by ${}'$ to
highest weight vectors of weight $\sharp\La$ 
of $V_q(\sharp\LPP)\otimes V_q(\sharp\LP)$.

\end{appendix}

\end{document}